\documentclass[12pt]{amsart}
\usepackage{amsmath,amssymb}
\usepackage{graphicx,multirow,rotating}

\numberwithin{equation}{section}
\newtheorem{thm}{Theorem}[section]
\newtheorem{lem}{Lemma}[section]
\newtheorem{prop}{Proposition}[section]
\newtheorem{cor}{Corollary}[section]

\newtheorem{rem}{Remark}[section]

\def\1{{{\mbox{${\rm{1\negthinspace\negthinspace I}}$}}}}
\newcommand{\eref}[1]{(\ref{#1})}

\newcounter{hypc}
\newcommand{\hypothese}[1]{\stepcounter{hypc}
\tag{$\mathbf{A_{\thehypc}}$}
\label{#1}}

\newcounter{condc}

\newcommand{\condition}[1]{\stepcounter{condc}
\tag{$\mathbf{R_{\thecondc}}$}
\label{#1}}

\begin{document}
\title[Nonparametric estimation in an error-in-variables model]{Nonparametric estimation of the
regression function in an errors-in-variables model}
\author[F. Comte]{F. Comte$^{*,1}$}\thanks{$^1$ Universit\'e Paris V,  MAP5, UMR CNRS
8145, 45  rue des Saints-P\`eres, 75 270 PARIS cedex 06, France.
email: fabienne.comte@univ-paris5.fr.}
\author[M.-L. Taupin]{M.-L. Taupin$^2$}\thanks{$^2$ IUT de Paris V et
Universit\'e d'Orsay, Laboratoire de Probabilit\'es, Statistique
et Mod\'elisation, UMR 8628, B\^{a}timent 425,91405 Orsay Cedex, France.- email:
Marie-Luce.Taupin@math.u-psud.fr}

\begin{abstract} We consider the regression model with errors-in-variables
where we observe $n$ i.i.d. copies of $(Y,Z)$ satisfying $Y=f(X)+\xi, \; Z=X+\sigma\varepsilon$, involving
independent and unobserved random variables $X,\xi,\varepsilon$.
The density $g$ of $X$ is unknown, whereas the density of
$\sigma\varepsilon$ is completely known. Using the observations $(Y_i,
Z_i)$, $i=1,\cdots,n$, we propose an estimator of the regression function
  $f$, built as the ratio of two penalized minimum contrast estimators
of $\ell=fg$ and $g$, without any prior knowledge on their smoothness. We prove that its $\mathbb{L}_2$-risk on a compact set is bounded by the sum of the two 
$\mathbb{L}_2(\mathbb{R})$-risks of   the estimators of $\ell$ and $g$, and
give the rate of convergence of such estimators for various smoothness classes
for $\ell$ and $g$, when the errors $\varepsilon$ are either ordinary smooth
or super smooth.
The resulting rate is optimal in a minimax sense in all cases where lower
  bounds are available.
\end{abstract}
\maketitle
\begin{center}
October 2004. Revised September 2005.
\end{center}

Errors-in-variables. Density deconvolution. Minimax estimation. Nonparametric regression.
Projection estimators.

\section{Introduction}

We consider that we observe  $n$ independent and  identically distributed
(i.i.d.) copies of $(Y,Z)$ satisfying the following
errors-in-variables regression model
\begin{eqnarray}\label{model}
\left\{\begin{array}{l}Y=f(X)+\xi \\ Z=X+\sigma\varepsilon,
\end{array}\right. \end{eqnarray}
involving independent and unobserved, random variables
$X,\xi,\varepsilon$ and an unknown regression function $f$. The
unobserved $X_i$'s, have common unknown density denoted by $g$.
The errors $\varepsilon_i$'s have common known density
$f_{\varepsilon}$, and $\sigma$ is the known noise level.
We assume moreover that all random variables have finite variance.
Our aim is to estimate the regression function $f$ on a compact
set denoted by $A$, by using the observations $(Y_i,Z_i)$ for
$i=1, \dots, n$, without any prior knowledge, neither on the
smoothness of $f$ nor on the smoothness of the density $g$. In
nonparametric errors-in-variables regression models, two factors
determine the estimation accuracy of the regression function:
first, the smoothness of the function $f$ to be estimated, and
second the smoothness of the errors density $f_\varepsilon$. As in the deconvolution framework, the worst rates of convergence
are obtained for the smoother errors density $f_\varepsilon$. In
this context, two classes of errors are considered: first the so
called ordinary smooth errors with polynomial decay of their
Fourier transform and second, the super smooth errors with Fourier
transform having an exponential decay.

Many papers deal with parametric or semi-parametric estimation in
errors in variables models, but we only mention here previous
known results in the general nonparametric case. In this context
most of the proposed estimators are some Nadaraya-Watson kernel
type estimators, constructed as the ratio of two deconvolution
kernel type estimators, see e.g. Fan \textit{et al.}~(1991), Fan
and Masry~(1992), Fan and Truong~(1993), Masry~(1993),
Truong~(1991), Ioannides and Alevizos~(1997).
One assumption usually done in all those works, is that the
regularity of the regression function $f$ and the regularity of
the density $g$ of the design are equal. In particular, when the
regression function $f$ and the density $g$ admit $k$th-order
derivatives, Fan and Truong~(1993) give upper and lower bounds of
the minimax risk for quadratic pointwise risk and for
$\mathbb{L}_p$-risk on compact sets for ordinary and super smooth
errors $\varepsilon$.


In a slightly different way, Koo and Lee~(1998) propose an
estimation method based on $B$-spline, when the errors are
ordinary smooth.  This method also relates to estimation of the
regression function as a ratio of two estimators.

To our knowledge, all previous papers consider that the regression
function and the density $g$  belong to the same smoothness class
and that this common class is known.

We propose here an estimation procedure of $f$, that does not
require any prior knowledge on the regularity of the unknown
functions $f$ and $g$.  Our estimation procedure is based  on the
classical idea that the regression function $f$ at point $x$ can
be written as the ratio
$$f(x)=\mathbb{E}(Y|X=x)=\frac{\int y f_{X,Y}(x,y)dy}{g(x)}=\frac{(fg)(x)}{g(x)},$$
with  $f_{X,Y}$ the joint density of $(X,Y)$. Hence $f$ is
estimated by a ratio of an adaptive estimator $\tilde{\ell}$ of
$\ell=fg$ and of an adaptive estimator $\tilde{g}$ of $g$, both of
them being built by minimization of penalized contrast functions.
The contrasts are determined by projection methods and the
penalizations give an automatic choice of the relevant projection
spaces.

We give upper bounds on the $\mathbb{L}_2$-risk on a compact set
for the regression function $f$ as well as for the
$\mathbb{L}_2(\mathbb{R})$-risk of the density $g$ when the errors
are either ordinary or super smooth. We show in particular that
the $\mathbb{L}_2$-risk on a compact set of our estimator $\tilde
f$ of $f$ is bounded by the sum of the risks
of $\tilde \ell$ and $\tilde g$. The rate of convergence of $\tilde f$ is thus
given by the slower rate between the rate of the adaptive
estimation of $g$ and the rate of the adaptive estimation of
$\ell=fg$. The resulting estimator automatically reaches the
minimax rates in standard cases where lower bounds are
available. The other cases are intensively discussed. In other
words, our procedure provides an adaptive estimator, in the sense
that its construction does not require any prior knowledge on the
smoothness of $f$ nor $g$, which seems often optimal.

The paper is organized as follows. In Section 2, we describe the
estimators. Section 3 is devoted to the presentation of the upper bounds for the resulting
$\mathbb{L}_2$-risks with some discussions about the optimality in the minimax sense of the estimators.
All proofs and technical lemmas are gathered in Section
4.

\section{Description of the estimators}
For $u$ and $v$ in $\mathbb{L}_2(\mathbb{R})$,
$u^*$ is the Fourier transform of $u$ with $u^*(x)=\int e^{itx}u(t)dt$,
$u*v$ is the convolution product, $u*v(x)=\int
u(y)v(x-y)dy$, and $ <u,v>=\int
u(x)\overline{v}(x)dx$ with $z\overline{z}=|z|^2.$ The quantities
$\|u\|_1$, $\|u\|_2$,  $\|u\|_\infty$ and 
$\|u\|_{\infty, K}$ denote
$\|u\|_1=\int |u(x)|dx,~ \|u\|_2^2=\int |u(x)|^2 dx,~ \|u\|_\infty=\sup_{x \in \mathbb{R}}|u(x)|,$
$\|u\|_{\infty, K}=\sup_{x \in K}|u(x)|$.

Subsequently we assume that $f_{\varepsilon}\in
\mathbb{L}_2(\mathbb{R})$, $f^*_{\varepsilon} \in
\mathbb{L}_2(\mathbb{R})$ with $f^*_{\varepsilon}(x)\not=0$ for
all $x\in\mathbb{R}$.

\subsection{Projection spaces}

Consider $\varphi(x)=\sin(\pi x)/(\pi x), \mbox{ and } \varphi_{m,j}(x) =
\sqrt{D_m} \varphi(D_mx-j).$
Here, we take $D_m=m$ and $m\in
\mathcal{M}_n=\{1,\cdots,m_n\}$, but when $D_m=2^m$, the basis  $\{\varphi_{m,j}\}_{j
\in \mathbb{Z}}$ is known as the Shannon basis. It is well known (see for instance Meyer (1990),
p.22),  that $\{\varphi_{m,j}\}_{j
\in \mathbb{Z}}$ is an orthonormal basis 
of the space $S_m$ of square
integrable functions having a Fourier transform with compact
support contained in $[-\pi D_m, \pi D_m]$ , that is $$S_m= {\rm Vect}\{\varphi_{_{m,j}}, \;
j\in \mathbb{Z}\}=\{f\in \mathbb{L}_2(\mathbb{R}), \mbox{ with }
\mbox{supp}(f^*) \mbox{ contained in }[-\pi D_m,\pi D_m ]\}.$$

Since the orthogonal projection of $g$ and $\ell$ on $S_m$, $g_m$
and $\ell_m$, $g_{m}=\sum_{j\in {\mathbb Z}} a_{m,j}(g)
\varphi_{m,j}$ and $\ell_{m}=\sum_{j\in {\mathbb Z}} a_{m,j}(\ell)
\varphi_{m,j}$ with $a_{m,j}(g) = <\varphi_{m,j},g>,$ and
$a_{m,j}(\ell) = <\varphi_{m,j},\ell>$, involve infinite sums, we
consider in practice, the truncated spaces $S_{m}^{(n)}$ defined
as
$$S_{m}^{(n)}= {\rm Vect }\left\{\varphi_{m,j}, |j|\leq k_n\right\}$$ where
$k_n$ is an integer to be chosen later. The family
$\{\varphi_{m,j}\}_{\vert j\vert \leq k_n}$ is an orthonormal
basis of $S_m^{(n)}$, and the orthogonal projection of $g$ and
$\ell$ on $S_{m}^{(n)}$ denoted by $g_{m}^{(n)}$ and
$\ell_{m}^{(n)}$, are given by $g_{m}^{(n)}=\sum_{|j|\leq k_n}
a_{m,j}(g) \varphi_{m,j}$ and $\ell_{m}^{(n)}=\sum_{|j|\leq k_n}
a_{m,j}(\ell) \varphi_{m,j}.$

\subsection{Construction of the minimum contrast estimators}

For
$r\in \mathbb{R}$ and $d>0$, we denote by $r^{(d)}=\mbox{sign}(r)\min(\vert
r\vert, d)$, and thus define the trimmed estimator of
$f$ by
\begin{eqnarray}\label{estimf} \hat f_{\breve m_\ell,\breve m_g}=(\hat
\ell_{\breve m_\ell}/\hat g_{\breve m_g})^{(a_n)},
\end{eqnarray}
with $a_n$ being suitably chosen, $\breve m_\ell$ and $\breve  m_g$ minimizing
the $\mathbb{L}_2(\mathbb{R})$ risks of $\hat \ell_{\breve m_\ell}$ the projection estimator on a space
$S_{\breve m_\ell}^{(n)}$, and of $\hat g_{\breve m_g}$ the projection estimator on a space
$S_{\breve m_g}^{(n)}$, defined as follows.

\noindent The estimator of $\ell$, is defined by 
\begin{equation}\label{truncsanssel}  \hat \ell_m = \arg\min_{t\in
S_m^{(n)}} \gamma_{n,\ell}(t),\end{equation}
with $\gamma_{n,\ell}$ defined, for $t\in S_m^{(n)}$, by
\begin{equation}\label{criteres}
  \gamma_{n,\ell}(t)=\|t\|^2-2n^{-1}\sum_{i=1}^n(Y_iu_t^*(Z_i)) \mbox{ with } u_t(x)=(2\pi)^{-1}
t^*(-x)/f_{\varepsilon}^*(-x),
\end{equation}
that is
$\hat \ell_m =
\sum_{|j|\leq k_n} \hat a_{m,j}(\ell) \varphi_{m,j}$ with
$\hat a_{m,j}(\ell) = n^{-1}\sum_{i=1}^n
Y_i  u^*_{\varphi_{m,j}}(Z_i)$.

By using Parseval and inverse Fourier formulas, we get that
$$ {\mathbb E}(Y_1u_t^*(Z_1)) = {\mathbb
E}(f(X_1)u_t^*(Z_1)) = \langle u_t^**f_{\varepsilon}, fg\rangle =
\frac 1{2\pi}\langle f_{\varepsilon}^* t^*/f_\varepsilon^*,
(fg)^*\rangle = \frac 1{2\pi} \langle t^*,(fg)^*\rangle =
\langle t,\ell\rangle.$$ Therefore, we find that
$\mathbb{E}(\gamma_{n,\ell}(t))= \|t\|^2_2-2\langle \ell, t\rangle
= \|t-\ell\|^2_2 -\|\ell \|^2_2$ which is minimal when $t=\ell$.
This shows that $\gamma_{n,\ell}(t)$ suits well for the estimation
of $\ell=fg$.

By using the estimation procedure described in Comte\textit{ et al.}~(2005a),
the estimator of $g$ on
$S_m^{(n)}$ is defined by $\hat g_m= \sum_{|j|\leq k_n} \hat a_{m,j}(g) \varphi_{m,j} \;
\mbox{ with } \; \hat a_{m,j}(g)= n^{-1}\sum_{i=1}^n
u^*_{\varphi_{m,j}}(Z_i),$ that is
\begin{eqnarray}
\label{estimg}\hat g_m=\arg\min_{t\in S_m^{(n)}}\gamma_{n,g}(t)\end{eqnarray}
with $\gamma_{n,g}$ defined, for $t\in S_m^{(n)}$ by
$\gamma_{n,g}(t) = \|t\|^2_2 -2n^{-1} \sum_{i=1}^n
u_t^*(Z_i),$ with $u_t$ defined in \eref{criteres}.

\begin{rem} {\rm 
The use of $r^{(d)}$
avoids
the problems that may occur when $\hat g_{m_2}$ takes small values.}
\end{rem}

\subsection{Construction of the minimum penalized contrast estimators}
In order to construct the minimum penalized contrast estimators, and especially
to define the penalty functions, we need to precise the behavior of
$f_\varepsilon^*$, described as follows. We assume that, for all $x$ in ${\mathbb R}$,
\begin{align}
&\hypothese{TFfeps}
\kappa_0(x^2+1)^{-\alpha/2}
\exp\{-\beta\vert x\vert^\rho\}\leq \vert f_\varepsilon^*(x)\vert \leq \kappa_0^{\prime}(x^2+1)^{-\alpha/2}
\exp\{-\beta\vert x\vert^\rho\}.
\end{align}
Only the left-hand side of \eref{TFfeps} is required to define the penalty
function and for upper bounds. The right-hand side is needed when we consider lower
bounds and the question of optimality in a minimax sense.
When $\rho=0$, $\alpha$ has to be such that $\alpha>1/2$ .
When $\rho=0$ in
\eref{TFfeps}, the errors are  usually called ``ordinary smooth'' errors, and ``super smooth''
errors when  $\rho>0.$ The standard examples are the following:
Gaussian or Cauchy distributions are super smooth of order ($\alpha=0$, $\rho=2$) and
($\alpha=0$, $\rho=1$) respectively, and
the double exponential distribution is ordinary smooth
($\rho=0$) of order $\alpha=2$.

By convention, we set $\beta=0$ when $\rho=0$ and we assume that $\beta>0$ when $\rho>0$.
 In the same way, if $\sigma=0$, the $X_i$'s are directly observed without noise and we set $\beta=\alpha=\rho=0$.

Under the assumption \eref{TFfeps}, the regression function $f$ is
estimated by $\tilde f$ defined as
\begin{equation}\label{defdeftilde} \tilde f= ({\tilde \ell}/{\tilde g})^{(a_n)},\end{equation}
where $\tilde \ell $ is the adaptive estimator defined
by
\begin{equation}\label{trunctilde}
\tilde \ell =\hat \ell_{\hat m_{\ell}} \mbox{ with } \hat
m_{\ell}= \arg\min_{m\in {\mathcal M}_{n, \ell}}
\left[\gamma_{n,\ell}(\hat \ell_m) + \; {\rm pen}_{\ell}(m)\right],
\end{equation}
$\tilde g$ is the adaptive estimator defined
 as in
Comte\textit{ et al.}~(2005a), by
\begin{equation}\label{truncgtilde}
\tilde g=\hat g_{\hat m_g} \mbox{ with } \hat m_g= \arg\min_{m\in
{\mathcal M}_{n,g}} \left[\gamma_{n,g}(\hat g_m) + \; {\rm
pen}_g(m)\right], \end{equation} where ${\mathcal M}_{n,\ell}$ and ${\mathcal
M}_{n,g}$ are some restrictions of ${\mathcal M}_n$ given below,  and where pen$_\ell$ and pen$_g$ are
data driven penalty functions given by
\begin{eqnarray}\label{penal}
\quad{\rm pen}_{\ell}(m) = \kappa'(\lambda_1+\mu_2)
[1+ \hat m_2(Y)]\tilde\Gamma(m)/n,\quad {\rm pen}_g(m)=
\kappa(\lambda_1+\mu_1)
\tilde\Gamma(m)/n,
\end{eqnarray}
with
\begin{equation}\label{estey2}
\hat m_2(Y) = \frac 1n\sum_{i=1}^n Y_i^2, \quad \mbox{ and }
\quad\tilde{\Gamma}(m)=D_m^{2\alpha+\max(1-\rho, \min((1+\rho)/2, 1))}
\exp\{2\beta\sigma^\rho(\pi D_m)^{\rho}\}.\end{equation}
The constants $\lambda_1, \mu_1$ and $\mu_2$ are some known constants, only depending
on $f_\varepsilon$ and $\sigma$ (assumed to be known), to be defined later (see (\ref{lambda1}),
(\ref{mu1}) and (\ref{mu2})), and $\kappa$ and $\kappa'$ are some numerical
constants.


\begin{rem}
\label{penalites}{\rm
First note that the penalty functions in \eref{penal} have
the same form with different constants. More precisely, in both cases, the
penalties are of order $D_m^{2\alpha+1-\rho}\exp(2\beta \sigma^\rho(\pi D_m)^\rho)$
if $0\leq \rho\leq 1/3$, $D_m^{2\alpha+(1+\rho)/2}\exp(2\beta \sigma^\rho(\pi D_m)^\rho)$ if
$1/3\leq \rho \leq 1$ and of order $D_m^{2\alpha+1}\exp(2\beta \sigma^\rho(\pi D_m)^\rho)$ if $\rho\geq 1$.

Second, the constants involve $\kappa$ and $\kappa'$,  universal
numerical constants, as well as constants $\lambda_1$, $\mu_1$, $\mu_2$ related to the known
errors density $f_\varepsilon$. Any constant greater than any
well chosen constant also suits for theoretical results. In practice,
such constants are usually calibrated by some intensive simulation
studies. We refer to Comte \textit{et al.}~(2005a, 2005b) for
further details on penalty calibration  as well as for details on
the implementation of such estimators in density deconvolution
problems. }\end{rem}

\section{Rates of convergence and adaptivity}

\subsection{Assumptions}\label{secmodel}
We consider Model \eref{model} under (\ref{TFfeps}) and the following additional assumptions.
\begin{align}
 &\ell \in \mathbb{L}_2(\mathbb{R})\mbox{ and }\ell\in\mathcal{L}=\left\{\phi \mbox{ such that }
\int x^2\phi^2(x)dx \leq \kappa_{\mathcal{L}} <\infty\right\},\hypothese{lL2}
\\
& f  \in
\mathcal{F}_G=\{\phi \mbox{ such that }\sup_{x\in G }|\phi(x)| \leq
\kappa_{\infty,G}<\infty\},\mbox{ where }G\mbox{ is the support of
}g.\hypothese{fbornee}
\\
 & \hypothese{gL2}g \in \mathbb{L}_2(\mathbb{R})\mbox{ and } g\in \mathcal{G}=\{\phi, \mbox{ density, such that } \int x^2 \phi^2(x)dx < \kappa_{\mathcal {G}} <\infty\}.\notag\\
 &\hypothese{gbornee} \mbox{ There exist } g_0, g_1 \mbox{ positive constants such that for all } x \in A,  g_0\leq g(x) \leq g_1.
\end{align}

Note that we do not assume that $g$ is compactly supported
but only that $f$ is bounded on the support of $g$. It follows that if $g$ is
compactly supported then $f$ has to be bounded on a compact set. But if $g$
has $\mathbb{R}$ as support then the regression function has to be bounded on
$\mathbb{R}$. We estimate $f$ only on a compact set denoted by
$A$. Hence, the assumption \eref{gbornee} implies that
$A\subset G$ and therefore under \eref{fbornee} and \eref{gbornee}, $f$ is
bounded on $A$. The assumptions \eref{fbornee} and \eref{gL2} imply that
\eref{lL2} holds, with $\kappa_{\mathcal{L}} = \kappa_{\infty,
  G}^2\kappa_{\mathcal{G}}$.

 Classically, the slowest rate of
convergence for estimating $f$ and $g$ are obtained for super
smooth errors density. In particular, when $f_\varepsilon$ is the
Gaussian density the minimax rate of convergence obtained by Fan
and Truong~(1993) when $f$ and $g$ have the same H\"olderian type
regularity is of order a power of $\ln(n)$. Nevertheless, those
rates can be improved by some additional regularity conditions on
$f$ and  $g$ described  as follows.
\begin{align}
\condition{super} \mathcal{
S}_{a,r,B}(C_1)=\{ \psi \in \mathbb{L}_2(\mathbb{R})\; : \; \mbox{ such that
} \int_{-\infty}^{+\infty}
|\psi^*(x)|^2(x^2+1)^{a}\exp\{2B |x|^{r}\} dx\leq
C_1\},\end{align}
for $a,r,B,C_1$ some nonnegative real numbers.
The smoothness class in \eref{super} is classically considered in
nonparametric estimation, especially in deconvolution.
When $r=0$, this corresponds to
Sobolev spaces of order $a$. The densities belonging to
$\mathcal{ S}_{a,r,B}(C_1)$ with $r>0, B>0$
are infinitely many times differentiable,
admit analytic continuation on a finite width strip when $r=1$
and on the whole complex plane if $r=2$.

\subsection{Risks bounds for the minimum contrast estimators}
We start by presenting some general bound for the risk.
\begin{prop}
\label{pestimgl} Consider the estimators $\hat \ell_{D_m}= \hat
\ell_m$ and $\hat g_{D_m}= \hat
g_m$ of $\ell$ and $g$ defined by \eref{truncsanssel} and \eref{estimg}.
Let $\Delta(m)=D_m \pi^{-1}\int_0^{\pi D_m}\vert
f_\varepsilon^*(D_mx\sigma)\vert^{-2}dx.$
Then, under \eref{lL2} and \eref{gL2},
\begin{equation}
\label{risktruncl1} {\mathbb E}(\|\ell-\hat \ell_m\|^2_2) \leq
\|\ell-\ell_m\|^2_2 + 2{\mathbb
E}(Y_1^2)\Delta(m)/n +(\kappa_{{\mathcal L}}+\parallel\ell\parallel_1)D_m^2/k_n
\end{equation}
and
\begin{equation}
\label{risktruncg} {\mathbb E}(\|g-\hat g_m\|^2_2) \leq
\|g-g_m\|^2_2 + 2\Delta(m)/n+(\kappa_{{\mathcal G}}+1)D_m^2/k_n.
\end{equation}
\end{prop}
As in deconvolution problems, the variance term
$\Delta(m)/n$ depends on the rate of decay of the Fourier
transform $f_\varepsilon^*$, with larger variance for fast decreasing
$f_\varepsilon^*$. Under \eref{TFfeps}, the variance term is bounded in the
following way
\begin{eqnarray}\label{gamma}
\Delta(m)\leq \lambda_1 \Gamma(m) \quad\mbox{where}\quad \Gamma(m)=D_m^{2\alpha+1-\rho}\exp(2\beta\sigma^{\rho}(\pi D_m)^{\rho}),
\end{eqnarray}
with 
\begin{eqnarray}
\label{lambda1} \quad\;\;\;\lambda_1= (\sigma^2\pi^2+1)^{\alpha}/(\pi^{\rho}
\kappa_0^{2}R(\beta,\sigma,\rho))\mbox{ with } R(\beta,\sigma, \rho)= \1_{\rho=0} + 2\beta \rho\sigma^{\rho} \1_{0<\rho\leq 1} + 2\beta \sigma^{\rho} \1_{\rho>1},
\end{eqnarray}
In order to ensure that $\Gamma(m)/n$ is bounded, we only consider models such that
$\pi D_m=m\leq m_n$ in  $\mathcal{M}_n=\{1,\cdots,m_n\}$ with
 \begin{eqnarray}
\label{mn}
m_n\leq \left\{
\begin{array}{ll}
\pi^{-1}n^{1/(2\alpha+1)} &\mbox{ if }\rho=0\\
\displaystyle\pi^{-1}\left[\frac{\ln(n)}{2\beta\sigma^\rho}+\frac{2\alpha+1-\rho}{2\rho\beta
  \sigma^\rho}\ln\left(\frac{\ln(n)}{2\beta\sigma^\rho}\right)\right]^{1/ \rho} &\mbox{ if
} \rho>0.
\end{array}
\right.
\end{eqnarray}

Lastly, the bias terms $\|\ell-\ell_m\|_2^2$ and $\|g-g_m\|^2_2$ depend, as usual, on
the smoothness of the functions $\ell$ and $g$. They have the expected order for
classical smoothness classes since they relate to the distance
between $g$ and the classes of entire functions having Fourier
transform compactly supported on $[-\pi D_m, \pi D_m]$ (see
Ibragimov and Hasminskii~(1983)).

Since $\ell_m$ and $g_m$ are  the orthogonal projections of
$\ell$ and $g$ on $S_m$, when $\ell$ belongs
$\mathcal{ S}_{a_\ell,r_\ell,B_\ell}(\kappa_{a_\ell})$ and
$g$ belongs
$\mathcal{ S}_{a_g,r_g,B_g}(\kappa_{a_g})$
 defined by \eref{super}, then
 \begin{equation}\label{biaiscarre}
\|\ell-\ell_m\|^2_2
= (2\pi)^{-1} \int_{|x|\geq \pi D_m} |\ell^*|^2(x)dx\leq [\kappa_{a_\ell}/(2\pi)]
(D_m^2\pi^2+1)^{-a_\ell}\exp\{-2B_\ell \pi^{r_\ell}D_m^{r_\ell}\},
\end{equation}
and the same holds for $\parallel g_m-g\parallel_2^2$ with
$(a_\ell,B_\ell,r_\ell)$ replaced by $(a_g,B_g,r_g)$.

\begin{cor}
Under \eref{TFfeps}, \eref{lL2} and \eref{gL2},  let $\Gamma(m)$ and $\lambda_1$
 being defined
in (\ref{gamma}) and (\ref{lambda1}).
Assume that
$k_n\geq n$, that $\ell$ belongs to $\mathcal{
S}_{a_\ell,r_\ell,B_\ell}(\kappa_{a_\ell})$ and that $g$ belongs
to $\mathcal{ S}_{a_g,r_g,B_g}(\kappa_{a_g})$
 defined by
\eref{super}. Then
$$
{\mathbb E}(\|\ell-\hat \ell_m\|^2_2)\leq \frac{\kappa_{a_\ell}}{2\pi} (D_m^2\pi^2+1)^{-a_\ell} e^{-2B_\ell \pi^{r_\ell}D_m^{r_\ell} } +
2\lambda_1{\mathbb
E}(Y_1^2)\Gamma(m)/n+D_m^2(\kappa_{\mathcal{L}}+\parallel\ell\parallel_1)/n
,$$
 and
$$
{\mathbb E}(\|g-\hat g_m\|^2_2)\leq\frac{\kappa_{a_g}}{2\pi} (D_m^2\pi^2+1)^{-a_g} e^{-2B_g
  \pi^{r_g}D_m^{r_g} } +
2\lambda_1\Gamma(m)/n+(\kappa_{\mathcal{G}}+1)D_m^2/n
.$$
\end{cor}

\begin{rem}
\label{choixkn} {\rm We point out that the $\{\varphi_{m,j}\}$ are
${\mathbb R}$-supported (and not compactly supported) and hence,
we obtain estimations of $\ell$ and $g$ on the whole line and not
only on a compact set as for usual projection estimators. This is
a great advantage of this basis even if, due to the truncation
$\vert j\vert \leq k_n$,  it induces the residual terms
$D_m^2(\kappa_{\mathcal{L}}+\parallel \ell\parallel_1)/k_n$ and
$D_m^2(\kappa_{\mathcal{G}}+1)/k_n$, in the upper bounds of the
risks. The most important thing is that the choice of $k_n$ does
not influence the other terms. Consequently, we can find a
relevant choice of $k_n$ ($k_n\geq n$ under \eref{lL2} and
\eref{gL2}), that makes those additional terms unconditionally
negligible with respect to the bias and variance terms. The
condition $k_n\geq   n$ allows us to
 construct truncated spaces $S_{m}^{(n)}$ using $O(n)$ basis vectors and hence
to use a tractable and fast algorithm. 
The choice of larger $k_n$, independent of $\ell$ and $g$, does not change the efficiency of our
estimator from a statistical point of view but will only change the speed of the
algorithm from a practical point of view.
}
\end{rem}




\begin{center}
{\small
\begin{table}[!ht]
\begin{tabular}{|c|l||c|c||}\cline{3-4}\cline{3-4}
\multicolumn{2}{c||}{} &\multicolumn{2}{c||}{$f_\varepsilon$} \\\cline{3-4}
\multicolumn{2}{c||}{} & $\rho=0$ & $ \rho>0$ \\
\multicolumn{2}{c||}{} & ordinary smooth & super smooth \\\hline\hline
\multirow{4}{.2cm}{\\\vfill\null $g$} & $\;$ & $\;$ & $\;$ \\
& $\begin{array}{l}
  r_\ell=0\\
  \small{\mbox{Sobolev}(s)}
\end{array}$ &
$\begin{array}{l}
  \pi D_{\breve m_\ell}=O(n^{1/(2\alpha+2a_\ell +1)})\\
  \mbox{rate}=O(n^{-2a_\ell/(2\alpha+2a_\ell+1)})
\end{array}$  &
$\begin{array}{l}
  \pi D_{\breve m_\ell}=[\ln(n)/(2\beta\sigma^{\rho}+1)]^{1/\rho}\\
  \mbox{rate}=O( (\ln(n))^{-2a_\ell/\rho})
\end{array}$  \\
$\;$ & $\;$ & $\;$ & $\;$ \\
\cline{2-4}
& $\begin{array}{l}
  r_\ell>0\\
  \mathcal{C}^\infty
\end{array}$ &
$\begin{array}{l} \\
  \pi D_{\breve m_\ell}=\left[{\ln(n)/2B_\ell}\right]^{1/r_\ell} \\
  \mbox{ rate}= \displaystyle  O\left(\frac{\ln(n)^{(2\alpha+1)/r_\ell}}n\right)\\
  \;\; \end{array}$ &
$\begin{array}{c}
  \pi D_{\breve m_\ell}  \mbox{ implicit solution of } \\
  {D_{\breve m_\ell}}^{\!\!\!\!2\alpha+2a_\ell+1-r_\ell}e^{2\beta \sigma^\rho
 (\pi D_{\breve m_\ell})^\rho+2B (\pi D_{\breve m_\ell})^{r_\ell}}\\
  \qquad= O(n)\\
\end{array}$
\\
\hline\hline
\end{tabular}\\~\\
 \caption{Best choices of $D_{\breve m_\ell}$ minimizing ${\mathbb E}(\|\ell-\hat
\ell_m\|^2_2)$ and resulting rates for $\hat \ell_{\breve m_\ell}$.}\label{rates}
\end{table}}
\end{center}
\vspace{-1cm}


%

For the case $r_\ell>0$ and $\rho>0$, the choice $\pi D_{\breve m_\ell}=[\ln(n)/(2\beta\sigma^\rho+1)]^{1/\rho}$ leads to a rate which is faster
than any power of $\ln(n)$ and slower than any power of $n$. For instance if $r_\ell=\rho$, the rate
is of order $[\ln(n)]^b n^{-B_\ell/(B_\ell+\beta\sigma^\rho)}$ with $b=[-2a_\ell \beta\sigma^\rho+(2\alpha -r_\ell+1)B_\ell]/[r_\ell
(\beta\sigma^\rho+B_\ell)]$.

The same table holds for $g$, by replacing $(a_\ell,B_\ell,r_\ell)$ by $(a_g,B_g,r_g)$.
For $D_{\breve m_g}$ chosen in the same way as $D_{\breve m_{\ell}}$ in Table 1, the rate of
convergence of $\hat g_{\breve m_g}$ is the minimax rate of convergence, as given in Fan~(1991a) for
$r_g=0$, in Butucea~(2004) for $r_g>0$ and $\rho=0$ and in Butucea and Tsybakov~(2004) for
$0<r_g<\rho$ and $a_g=0$.




The rate of convergence of $\hat f_{\breve m_\ell,\breve m_g}$ is given by the
following proposition.
\begin{prop}
\label{vitf}
Under \eref{TFfeps}, \eref{lL2}, \eref{fbornee}, \eref{gL2}, and \eref{gbornee}, assume
that $g$ belongs to some space $\mathcal{ S}_{a_g,r_g,B_g}(\kappa_{a_g})$ defined by \eref{super}
with $a_g>1/2$ if $r_g=0$.  Let $\hat f_{\breve m_\ell,\breve m_g}$ be defined
by \eref{estimf}, with $\breve m_\ell$ and $\breve m_g$ such that $D_{\breve
  m_{\ell}}$ and $D_{\breve m_g}$
minimize the risks ${\mathbb E}(\|\ell-\hat \ell_{m} \|_2^2)$ and ${\mathbb E}(\|g- \hat
g_{m}\|_2^2)$ respectively. If $a_n=n^k$ for $k>0$, and $k_n\geq
n^{3/2}$, then, for $n$ great enough and $C_0=Kg_0^{-2}(1+g_1g_0^{-2}\kappa_{\infty,G})$,
\begin{equation}
\label{majhatf}
\hspace{1cm}\mathbb{E}\|(\hat f_{\breve m_\ell,\breve m_g}-f)\1_{A}\|_2^2\leq
C_0[{\mathbb
E}(\|\ell-\hat \ell_{\breve m_\ell} \|_2^2)+{\mathbb E}(\|g- \hat
g_{\breve m_g}\|_2^2)]+o(n^{-1}).
\end{equation}
\end{prop}

\noindent
If $a_g\leq 1/2$ then  we only have a result of type
$\|(f-\hat f_{\breve m_\ell,\breve m_g})\1_A\|^2_2=O_p(\|\ell-\hat\ell_{\breve
  m_\ell}\|^2_2+\|g-\hat g_{\breve m_g}\|^2_2).$
Also note that the result holds when the constant $\kappa_{\infty,G}$ is replaced by $\parallel f\parallel_{\infty,A}$ if $f$ is bounded on the compact set $A$.

The performance of $\hat f_{\breve m_\ell,\breve m_g}$ is given by
the worst performance between the one of $\hat\ell_{\breve m_\ell}$ and
the one of $\hat g_{\breve m_g}$. Let us be more precise in some
examples. Under  the assumptions of Proposition \ref{vitf}:
\begin{itemize} \item If the $\varepsilon_i$'s are ordinary smooth,
\begin{itemize}
\item If $r_{\ell}=r_g=0$ and $\pi D_{\breve m_\ell}=O\left(
    n^{1/(2a_{\ell}+2\alpha+1)}\right)$ and $\pi D_{\breve m_g}=O\left( n^{1/(2a_g+2\alpha+1)}\right)$, then
$${\mathbb E}(\|(f-\hat f_{\breve m_\ell,\breve m_g})\1_A\|^2_2) \leq O(n^{-2a^*/(2a^*+2\alpha+1)}) \;\; \mbox{ with } \;\;
a^*= \inf(a_{\ell}, a_g).$$
\item If $r_{\ell}>0$, $r_g>0$, $\pi D_{\breve
    m_\ell}=(\ln(n)/2B)^{1/r_{\ell}}$ and $\pi D_{ \breve m_g}=(\ln(n)/2B)^{1/r_g}$, then
$${\mathbb E}(\|(f-\hat f_{\breve m_\ell,\breve m_g})\1_A\|^2_2) \leq O\left(
\frac{\ln(n)^{(2\alpha+1)/r^*}} n \right) \;\; \mbox{ with } \;\;
r^*= \inf(r_{\ell}, r_{g}).$$
\end{itemize}
\item If the $\varepsilon_i$'s are super smooth and $r_\ell=r_g=0$, $\pi
D_{ \breve m_\ell}=\pi D_{\breve m_g}=[\ln(n)/(2\beta\sigma^\rho+1)]^{1/\rho}$, then
$${\mathbb E}(\|(f-\hat f_{\breve m_\ell,\breve m_g})\1_A\|^2_2) \leq O( [\ln(n)]^{-2a^*/\rho})
\;\; \mbox{ with } \;\; a^*= \inf(a_{\ell}, a_g).$$
\end{itemize}

Since $\ell=fg$, the smoothness properties of $\ell$ are related to those of $f$ and of $g$.

When $\ell$ belongs to
${\mathcal S}_{a_\ell,0,B_\ell}(\kappa_{a_\ell})$ and $g$ belongs
to ${\mathcal S}_{a_g,0,B_g}(\kappa_{a_g})$ with $a_{\ell}=a_g$, then the
resulting rate is the
minimax rate given in Fan and Truong~(1993)
for H\"olderian regression functions and densities with the same regularity. It follows that
our estimator seems then optimal  in that case.
It is easy to see that the estimator is also optimal if
$a_g\geq a_\ell$, that is when the density $g$ is
smoother than the regression function $f$. But the optimality of
the rate of $\hat{f}_{\breve m_\ell,\breve m_g}$ when $a_\ell >a_g$,
that is when the regression function $f$ is smoother than $g$,
remains an open question. This is a known drawback of
Nadaraya-Watson type estimators for regression functions,
constructed as ratio of estimators. In ``classical''
regression models, when the $X_i$'s are observed, a lot of
methods, like local polynomial estimators, mean square
estimators..., avoid the need of regularity conditions on $g$ for
the estimation of $f$. The point is that standard methods solving
the regression problem do not seem to work in the errors-in-variables
model and it is an open problem to build an estimator of $f$ that does not
require the estimation of the density $g$.

>From the above results we see that  the choice of the
dimensions $D_{ \breve m_\ell}$ and $D_{\breve m_g}$ that realize the best
trade-off between the squared bias and the variance terms depends
on the unknown regularity coefficients of the functions $\ell$ and
$g$. In the next section we provide the upper bounds of the  risks
of the
penalized estimators, constructed without such smoothness knowledge.

\subsection{Risks bounds of the minimum penalized contrast estimators: adaptation}
\label{adaptation}

\begin{thm}\label{adaptatifgl}
Under the assumptions \eref{TFfeps}, \eref{lL2} and \eref{gL2}, let
\begin{equation}\label{mu1} \mu_1=\left\lbrace
\begin{array}{ll} 0 & \mbox{ if } \rho<1/3 \\
\beta(\sigma \pi)^{\rho}\lambda_1^{1/2}(\alpha,\kappa_0,\beta,\sigma,\rho)(1+\sigma^2\pi^2)^{\alpha/2}
\kappa_0^{-1}(2\pi)^{-1/2} &\mbox{ if } 1/3\leq\rho\leq
1,\\
\beta(\sigma \pi)^{\rho}\lambda_1(\alpha,\kappa_0,\beta,\sigma,\rho) &\mbox{ if }
\rho >1.\end{array}\right.\end{equation}
and
\begin{equation}\label{mu2}
\mu_2=
\mu_1 \1_{\{0\leq \rho <1/3\} \cup \{\rho>1\} } + \mu_1\|f_{\varepsilon}\|_2 \1_{\{1/3 \leq \rho\leq
1\}}.
\end{equation}
Let $k_n \geq
n$,  $\tilde \ell=\hat \ell_{\hat m_\ell}$  and $\tilde g = \hat
g_{\hat m_g}$ be defined by \eref{trunctilde} and (\ref{truncgtilde}) and  with {\rm pen}$_\ell$ and pen$_g$ given by
(\ref{penal}), for $\kappa$  and $\kappa'$ two universal numerical constants and
$1\leq m\leq m_n$, $m_n$ satisfying \eref{mn} and, if $\rho >0$, 
\begin{eqnarray}
\label{mn+}
m_n\leq
\pi^{-1}\left[\frac{\ln(n)}{2\beta\sigma^\rho}+\frac{2\alpha+\min[(1/2+\rho/2),1]}{2\rho\beta
  \sigma^\rho}\ln\left(\frac{\ln(n)}{2\beta\sigma^\rho}\right)\right]^{1/ \rho}.
\end{eqnarray}
\textbf{1)} \textbf{Adaptive estimation of $\mathbf{g}.$} (Comte\textit{ et
  al.}~(2005a)).\\
Then  $\tilde g$ satisfies ${\mathbb E}(\|g- \tilde g\|^2_2) \leq K\inf_{m\in {\mathcal M}_{n,g}}
\left[\|g-g_m\|^2_2 + D_m^2(\kappa_{\mathcal{G}}+1)/n+{\rm pen}_g(m)\right] + c/n$  where $K$ is a constant
and $c$ is a constant depending on
$f_{\varepsilon}$ and $A_g$.

\noindent\textbf{2)} \textbf{Adaptive estimation of $\mathbf{\ell}$.} Under
the assumption
\eref{fbornee}, if $\mathbb{E}\vert \xi_1\vert^8<\infty$
then $\tilde \ell$ satisfies $${\mathbb E}(\|\ell-\tilde \ell\|^2_2) \leq
K'\inf_{m\in {\mathcal M}_{n,\ell}} \left[\|\ell-\ell_m\|^2_2
+ D_m^2(\kappa_{\mathcal{L}}+\parallel \ell\parallel_1)/n
+ {\mathbb E}({\rm pen}_{\ell}(m))\right] +c'/n$$ where $K'$ is a constant and $c'$ is a constant depending on
$f_{\varepsilon}$, $\kappa_{\mathcal{L}}$, and $\|\ell\|_1$.
\end{thm}

\begin{rem}
{\rm
In Theorem \ref{adaptatifgl}, the penalty is random since it involves the term
$\hat m_2(Y)$, instead of the unknown quantity ${\mathbb E}(Y_1^2)$ which appears first.
The only price to pay for this substitution is the moment
condition $\mathbb{E}\vert \xi_1\vert^8<\infty$ instead of $\mathbb{E}\vert \xi_1\vert^6<\infty$ if
${\mathbb E}(Y_1^2)$ was in the penalty. Moreover, the term ${\mathbb E}({\rm pen}_{\ell}(m))$ in the bound is equal to pen$_{\ell}(m)$ with $\hat m_2(Y)$ replaced by ${\mathbb E}(Y_1^2)$. 
}\end{rem}

\begin{rem}
\label{penalites2}
{\rm According to Remark \ref{penalites}, the penalty functions are of order $\Gamma(m)/n$
if $0\leq \rho\leq 1/3$, of order $D_m^{3\rho/2-1/2}\Gamma(m)/n$ if
$1/3\leq \rho \leq 1$ and of order $D_m^{\rho}\Gamma(m)/n$ if $\rho\geq 1$.
When $\rho>1/3$, the penalty functions pen$_\ell(m)$ and pen$_\ell(m)$ have not exactly the order of the
variance $\Gamma(m)/n$, but a loss of order $D_m^{\min[(3\rho/2-1/2),\rho]}$
occurs, that is of order $D_m^{(3\rho-1)/2}$ if $1/3<\rho\leq
1$ and of order $D_m^{\rho}$ if $\rho>1$.
}
\end{rem}

\begin{rem} \textbf{Rates of convergence of $\mathbf{\tilde g}$}.
{\rm The rate of convergence of $\tilde{g}$ is the rate of convergence
of $\hat{g}_ {\breve m_g}$ when $0\leq \rho\leq 1/3$ or when $\rho>1/3$ and $r_g=0$
or $r_g<\rho$. And there is a logarithmic loss, as a price to pay for
adaptation when $r_g\geq\rho>1/3$. We refer to Comte\textit{ et al.}~(2005a)
for further comments on the optimality in a minimax sense of $\tilde{g}$.}
\end{rem}
\begin{rem} \textbf{Rates of convergence of $\mathbf{\tilde \ell}$}.
{\rm The rates, similar to the rates of $\tilde{g}$, are easy to deduce
from Theorem \ref{adaptatifgl} as soon as $\ell=fg$ belongs to some
smoothness class, but the procedure can reach the rate of
$\hat{\ell}_{\breve m_\ell}$, that uses the unknown smoothness parameter. If pen$_{\ell}(m)$ has the same order as the
variance order $\Gamma(m)/n$, then Theorem \ref{adaptatifgl}
guarantees an automatic trade-off between the squared  bias term
$\|\ell-\ell_m\|^2_2$ and the variance term, up to some
multiplicative constant. Else, there is some loss due to the
adaptation. Let us be more precise.

If $0\leq \rho\leq 1/3$, the errors $\varepsilon_i$'s are ordinary smooth or
super smooth with $\rho\leq 1/3$. If $\ell$ satisfies \eref{super}, the squared bias is bounded by applying \eref{biaiscarre}
which combined with the value of pen$_\ell(m)$, of order $ \Gamma(m)/n$ (see \eref{gamma})
gives that the estimator $\tilde g$ automatically reaches the best rate achievable
by the estimator $\hat{\ell}_{\breve m_\ell}$, as given in Table 1.

If $\rho>1/3$ the penalty function pen$_\ell(m)$ is slightly bigger
than the variance order $\Gamma(m)/n$. 
The rate of convergence remains the best rate if the bias $\|\ell-\ell_m\|^2_2$ is the dominating
term in the trade-off between $\|\ell-\ell_m\|_2^2$ and $\mbox{pen}_{\ell}(m)$.
When $r_\ell=0$ and $\rho>0$, the rate of
order $(\ln(n))^{-2a_\ell/\rho}$ is given by the bias term, and
the loss in the penalty function does not change the rate of the adaptive estimator $\tilde \ell$, which remains the best achievable
rate $\mathbb{E}\parallel \ell - \hat \ell_{\breve m_\ell}\parallel_2^2$.
In the same way, when $0<r_\ell<\rho$, the rate
is given by the bias term and thus this loss does
not affect the rate of convergence of $\tilde \ell$ either.

Let us now focus our discussion on the case where $\mbox{pen}_{\ell}(m)$
can be the dominating term in the trade-off between $\|\ell-\ell_m\|^2_2$
and $\mbox{pen}_{\ell}(m)$, that is when $r_\ell\geq \rho>1/3$. In that
case, there is a loss of order $D_m^{\min[(3\rho/2-1/2),\rho]}$ in the penalty function, compared to the variance term.
But this happens in cases where the order of the optimal
$D_m$ is less than $(\ln{n})^{1/\rho}$ and consequently the loss
in the rate is at most of order $\ln{n}$, when the rate is faster
than logarithmic: therefore the loss appears only in cases where
it can be seen as negligible.

In particular, there is no price to pay for the adaptation if the $\xi_i$'s are Gaussian and the
$\varepsilon_i$'s are ordinary smooth. Indeed, in that case,
the rate of convergence of the penalized estimator $\tilde{\ell}$, without any knowledge on $\ell$ or $g$,
 is the same as the rate given by the non penalized estimator $\hat{\ell}_{
   \breve m_\ell}$, requiring the knowledge of smoothness parameters.
But, if both the $\xi_i$'s and
the $\varepsilon_i$'s are Gaussian, then $\rho=2$
and a logarithmic negligible loss appears in the rate of $\tilde{\ell}$
compared to the rate of $\hat{\ell}_{\breve m_\ell}$.
}
\end{rem}

\begin{thm}\textbf{Adaptive estimation of $f$.}\label{adaptf}
Under the assumptions \eref{TFfeps}, \eref{lL2},  \eref{fbornee}, \eref{gL2}
and \eref{gbornee}, let $\tilde f$ be defined by (\ref{defdeftilde}) with $\tilde g$ and $\tilde \ell$ be
defined in \eref{truncgtilde} and \eref{trunctilde} with $\hat m_g\in
\mathcal{M}_{n,g}$ satisfying \eref{mn} and \eref{mn+},  $D_{m_{n,g}}\leq
(n/\ln(n))^{1/(2\alpha+2)}$ and $\hat m_{\ell} \in \mathcal{M}_{n,\ell}$ satisfying \eref{mn} and \eref{mn+}.
Assume that $g$ belongs to some space $\mathcal{
  S}_{a_g,r_g,B_g}(\kappa_{a_g})$ defined by \eref{super} with $a_g>1/2$ if $r_g=0$, and
that $\mathbb{E}\vert \xi_1\vert^8<\infty$.
If $k_n\geq n^{3/2}$, $a_n=n^k$ for $k>0$, for $n$ large enough,
$C_0=8Kg_0^{-2} $ and $C_1=4K'g_0^{-2}(2g_1^2+1)\kappa^2_{\infty,G}$, then
\begin{eqnarray}\nonumber {\mathbb E}(\|(f-\tilde f)\1_A\|^2_2)& \leq &
C_0\inf_{m\in {\mathcal M}_{n,\ell}}
[\|\ell-\ell_m\|^2_2 + D_m^2(\kappa_{\mathcal{L}}+\parallel \ell\parallel_1)/n+{\rm pen}_{\ell}(m)] \\
\label{doubleinf} && +
C_1\inf_{m\in {\mathcal M}_{n,g}} [\|g-g_m\|^2_2 + D_m^2(\kappa_{\mathcal{G}}+1)/n+{\rm pen}_g(m)] +  c/n
\end{eqnarray}
where $K$ and $K'$ are constants depending on $f_{\varepsilon}$, and $c$ is a constant depending on
$f_{\varepsilon}$, $f$ and $g$.
\end{thm}
\noindent As in Theorem \ref{estimf}, if $a_g\leq 1/2$ then it may happen that $D_{\tilde m_g}\geq
n^{1/(2\alpha+2)},$ and in this case we only have a result in probability:
$\|(f-\tilde f)\1_A\|^2_2=O_p(\|\ell-\tilde \ell\|^2_2+\|g-\tilde g\|^2_2).$
Moreover, the result holds when the constant $\kappa_{\infty,G}$ is replaced by $\parallel f\parallel_{\infty,A}$ if $f$ is bounded on the compact set $A$.
Also note that the remark \ref{choixkn} is still valid for all adaptive estimators.\\

\noindent {\bf Comments about the resulting rates for estimating $f$.}
First the rate of convergence of $\tilde{f}$ is given by the worst rate of
convergence between the rate of $\tilde{\ell}$ and $\tilde{g}$. Obviously all
the comments about $\hat{f}_{\breve m_\ell,\breve m_g}$, related to this fact keep holding here.

When $0\leq \rho<1/3$ or when $r_\ell\leq
\rho$ and $r_g\leq \rho$, then $\tilde{f}$ achieves the rate
of convergence of $\hat{f}_{\breve m_\ell, \breve m_g}$, given by the worst rate
of convergence between $\mathbb{E}\parallel
\hat{\ell}_{\breve m_\ell}-\ell\parallel_2^2$ and $\mathbb{E}\parallel
\hat{g}_{ \breve m_g}-g\parallel_2^2$. And when $r_g>\rho>1/3$ or
$r_\ell>\rho>1/3$, there is a logarithmic loss in the rate of
convergence of $\tilde{f}$ compared to the rate of convergence of
$\hat{f}_{\breve m_\ell,\breve m_g}$.

Since the regularity of $\ell$ is by definition the regularity of $fg$,  the rate of convergence of $\tilde{\ell}$ in fact depends on smoothness properties of $f$ and $g$.
As  a consequence, if $\ell$ and $g$  belong respectively to ${\mathcal
S}_{a_\ell,r_\ell,B_f}(\kappa_{a_\ell})$ and ${\mathcal
S}_{a_g,r_g,B_g}(\kappa_{a_g})$, then the rate of convergence of
$\tilde f$ is the rate of $\hat{f}_{\breve m_\ell,\breve m_g}$ when $0\leq \rho\leq 1/3$.
According to  Fan and Truong~(1993), this rate seems the minimax rate when
$a_\ell\leq a_g$ and $r_\ell=r_g=0$.
In the other cases, the question of the optimality in a minimax sense remains
open. Even if the regression function is smoother than $g$ and $0\leq \rho\leq
1/3$, the rate of convergence of $\tilde f$ has the order of the rate of
convergence of $\hat{f}_{\breve m_\ell,\breve m_g}$, but we do not know if the rate of
$\hat{f}_{\breve m_\ell,\breve m_g}$ is the minimax rate (see comments following Theorem \ref{estimf}).
When $\rho>1/3$, a loss appears between the rate of convergence of
$\tilde f$ and  the rate of convergence of $\hat{f}_{\breve m_\ell,\breve m_g}$. This loss
only appears, when $r_\ell>\rho$ or $r_g>\rho$ (see the comments after Theorem
\ref{adaptatifgl}), in cases where it is negligible with respect to the
rate.

\begin{rem}
\label{sigma} {\rm Obviously, the resulting rates for all
estimators depend on the noise level $\sigma$. The first point is
to note that if $\sigma=0$, then by convention
$B=0=\rho=0,\lambda=1$, and $Z=X$ is observed. In that case, $\Gamma(m)/n$ of
order $D_m/n$ has the
expected order for the variance term in ``usual regression'', when
the explanatory variables are observed, and the same holds for
the penalties ${\rm pen}_\ell$ and ${\rm pen}_g$. This
order $D_m/n$ is the expected penalty order for density estimation and nonparametric regression
estimation, when there is one model per dimension, as in our case.

The second point is to note that if $\sigma$ is small, then the
procedure automatically selects a dimension $D_m$ closed to the
dimension that would be selected in ``usual'' density estimation
and nonparametric regression estimation.}
\end{rem}

\section*{Concluding remarks}
Our estimation procedure provides an adaptive estimator in the sense that its
construction does not require any prior knowledge on the smoothness parameters
of the regression function $f$ and of the density $g$. This estimation procedure  allows to consider
various smoothness classes for the regression function and for the
density $g$ when the errors are either ordinary smooth or super smooth,
and to give upper bounds for the risk in all the cases.

The resulting rates of convergence for the estimation of $f$ are given by the
worst between the rate for the estimation of $fg$ and the rate for the
estimation of $g$. Nevertheless, they are the minimax
rates in cases where lower bounds are available. In the other cases, the
resulting rates are in most cases the best rates achievable if the smoothness
parameters were known. Some logarithmic loss, negligible compared to the
order of the rate, appears,  as a price to pay for the adaptation, when both the errors density
 $f_\varepsilon$ and $f g$ are super smooth with $f_\varepsilon$ strictly
 smoother than $fg$. This logarithmic loss appears when the influence of the noise $\sigma \varepsilon $ dominates the smoothness properties of $f$ and $g$.

\section{Proofs}\label{Proofs}
\subsection{Proof of Proposition \ref{pestimgl}}
By applying Definition (\ref{truncsanssel}), for any $m$
belonging to ${\mathcal M}_n$, $\hat \ell_m$ satisfies
$\gamma_{n,\ell}(\hat\ell_m)-\gamma_{n,\ell}(\ell_m^{(n)})\leq 0.$
Denoting by $\nu_n(t)$ the centered empirical process,
\begin{equation}\label{nun} \nu_n(t)=\frac 1n \sum_{i=1}^n \left(Y_iu_{t}^*(Z_i)-\langle
t, \ell \rangle\right),\end{equation}
and by using that $t\mapsto u_{t}^*$ is linear we get the following decomposition
\begin{eqnarray}
\label{difgamma}
\gamma_{n,\ell}(t)-\gamma_{n,\ell}(s)=\|t-\ell\|^2_2-\|s-\ell\|^2_2-2\nu_n(t-s)
\end{eqnarray}
and therefore, since by Pythagoras Theorem, $\|\ell -\ell_m^{(n)}
\|^2_2 = \|\ell - \ell_m\|^2_2 +\|\ell -\ell_m^{(n)}\|^2_2$, we infer that
$\|\ell-\hat \ell_m\|^2_2\leq  \|\ell-\ell_m\|^2_2
+ \|\ell_m-\ell_m^{(n)}\|_2^2+ 2\nu_n(\hat \ell_m -
\ell_m^{(n)})$. Using that $\hat
a_{m,j}(\ell)-a_{m,j}(\ell)= \nu_n(\varphi_{m,j})$, we get
\begin{equation}\label{nuphi}
\nu_n(\hat \ell_m-\ell_m^{(n)})=\sum_{|j|\leq k_n} (\hat
a_{m,j}(\ell)-a_{m,j}(\ell))\nu_n(\varphi_{m,j}) = \sum_{|j|\leq
k_n} [\nu_n(\varphi_{m,j})]^2,\end{equation} and consequently
\begin{eqnarray}
\label{decompnu2} {\mathbb E}\|\ell-\hat \ell_m\|^2_2\leq
\|\ell-\ell_m\|^2_2 +\|\ell_m-\ell_m^{(n)}\|^2_2+ 2\sum_{j\in
\mathbb{Z}}{\rm Var}[\nu_n(\varphi_{m,j})].\end{eqnarray}
 Now, since the  $(Y_i,Z_i)$'s are independent,
${\rm Var}[\nu_n(\varphi_{m,j})]=
n^{-1} {\rm Var}
[Y_1u_{\varphi_{m,j}}^*(Z_1)],$
and, arguing as in Comte\textit{ et al.}~(2005a), by  using
Parseval's formula we get that
\begin{eqnarray}
\label{var} \sum_{j\in{\mathbb Z}} {\rm Var}[\nu_n(\varphi_{m,j})]
\leq n^{-1}\parallel\sum_{j\in \mathbb{Z}}
|u^*_{\varphi_{m,j}}|^2\parallel_\infty { \mathbb E}(Y_1^2)\leq { \mathbb E}(Y_1^2)\Delta(m)/n.
\end{eqnarray}
 where
$\Delta$ is defined in Proposition \eref{pestimgl}.
Let us study the residual
term $\|\ell_m-\ell_m^{(n)}\|^2_2$, by simply writting that
$$\|\ell_m-\ell_m^{(n)}\|^2_2=\sum_{\vert j\vert \geq k_n} a_{m,j}^2(\ell)\leq
(\sup_{j}j a_{m,j}(\ell))^2\sum_{\vert j\vert \geq k_n }j^{-2}.$$
Now by definition
\begin{eqnarray*}
&&j a_{m,j}(\ell) =j\sqrt{D_m}\int \varphi(D_mx-j)\ell(x)dx\\
&\leq & D_m^{3/2}\int \vert x\vert \vert \varphi(D_mx -j)\vert
|\ell(x)|dx+\sqrt{D_m}\int \vert D_mx-j\vert \vert \varphi(D_mx
-j)\vert |\ell(x)|dx\\
\hspace{-1cm}&\leq & D_m^{3/2}\left( \int \vert \varphi(D_mx
-j)\vert^2 dx\right)^{1/2}\kappa_{\mathcal{L}}^{1/2}\!\!\!+\sqrt{D_m}\sup_{x}\vert
x\varphi(x)\vert \|\ell\|_1.
\end{eqnarray*}
Consequently
$j a_{m,j}\leq
D_m\|\varphi\|_2\kappa_{\mathcal{L}}^{1/2}+\sqrt{D_m}\|\ell\|_1/\pi,$ and
$\|\ell_m-\ell_m^{(n)}\|^2_2\leq \kappa
  (\kappa_{\mathcal{L}}+\|\ell\|_1^2) D_m^2/k_n$. $\hfill \Box$

\subsection{Proof of Proposition \ref{vitf}}
The proof of Proposition \ref{vitf} being rather similar to the proof of
Theorem \ref{adaptf} is omitted.
We refer to  Comte and Taupin~(2004) for further details.

\subsection{Proof of Theorem \protect \ref{adaptatifgl}}
We only prove the result with ${\mathbb E}(Y^2)$ in the
penalty instead of $\hat m_2(Y)$ and refer to Comte and Taupin~(2004) for the
complete proof with $\hat m_2(Y)$, as an application of Rosenthal's inequality
(see Rosenthal~(1970)).

For the study of $\tilde{\ell}$, the main difficulty compared to
the study of $\tilde{g}$ comes from the unbounded noise $\xi_i$. By definition, $\tilde \ell$
satisfies that for all $m\in {\mathcal M}_{n, \ell}$,
$\gamma_{n,\ell}(\tilde \ell)+ {\rm pen}_{\ell}(\hat m)\leq
\gamma_{n,\ell}(\ell_m^{(n)}) + {\rm pen}_{\ell}(m).$ Therefore, by
applying \eref{difgamma} we get that
\begin{eqnarray}
\parallel \tilde \ell-\ell\parallel_2^2 \leq \parallel\ell-\ell_m^{(n)}\parallel_2^2 +2 \nu_n(\tilde \ell -\ell_m^{(n)})
+ {\rm
pen}_{\ell}(m) -{\rm pen}_{\ell}(\hat m)
\label{basedec} .
\end{eqnarray}
Next, we use that if $t=t_1+t_2$ with $t_1$ in $S_m$ and $t_2$ in
$S_{m'}$, then $t$ is such that $t^*$ has its support in $[-\pi
D_{\max(m,m')}, \pi D_{\max(m,m')}]$ and therefore $t$ belongs to
$S_{m^*}$ where $m^*=\max(m,m')$. Denote by $B_{m, m'}(0,1)$ the set
$$B_{m, m'}(0,1)=\{t\in S_{\max(m,m')}^{(n)}
\;/\; \|t\|_2=1\}.$$ It follows that $$|\nu_n(\tilde \ell-\ell_m^{(n)})
|\leq \|\tilde \ell- \ell_m^{(n)}\|_2\sup_{t\in B_{m,\hat
m}(0,1)}|\nu_n(t)|,$$ where $\nu_n(t)$ is defined by (\ref{nun}).
Consequently, by using that
$2ab\leq x^{-1}a^2+xb^2$
\begin{eqnarray*} \|\tilde \ell-\ell\|_2^2 &\leq& \|\ell_m^{(n)}
-\ell\|_2^2 + \frac{1}{x}\|\tilde \ell-\ell_m^{(n)}\|_2^2  + x\sup_{t\in B_{m,\hat
m}(0,1)}\nu_n^2(t)+ {\rm pen}_{\ell}(m)- {\rm pen}_{\ell}(\hat m)
\end{eqnarray*}
and therefore, writing that $\|\tilde \ell-\ell_m^{(n)}\|_2^2\leq
(1+y^{-1})\|\tilde \ell-\ell\|_2^2+ (1+y)\|\ell-\ell_m^{(n)}\|_2^2$, with
$y=(x+1)/(x-1)$ for $x>1$, we infer that
\begin{eqnarray*} \|\tilde \ell-\ell\|_2^2
\leq \left(\frac{x+1}{x-1}\right)^2 \|\ell-\ell_m^{(n)}\|_2^2 +
\frac{x(x+1)}{x-1}\sup_{t\in B_{m,\hat m}(0,1)}\nu_n^2(t)
+\frac{x+1}{x-1} ({\rm pen}_{\ell}(m)- {\rm pen}_{\ell}(\hat m)).
\end{eqnarray*}
Choose some positive function $p_\ell(m,m')$ such that $x p_\ell(m,m')\leq
{\rm pen}_\ell(m) + {\rm pen}_\ell(m')$. Then, by denoting by $\kappa_x=(x+1)/(x-1)$,
\begin{eqnarray}\nonumber \|\tilde \ell-\ell\|_2^2
&\leq & \kappa_x^2 \|\ell-\ell_m^{(n)}\|_2^2+
x\kappa_x[\sup_{t\in B_{m,\hat
m}(0,1)}|\nu_n|^2(t)-p(m,\hat m)]_+ \\  &&
+\kappa_x\left(x p_{\ell}(m,\hat m) +{\rm pen}_{\ell}(m)- {\rm pen}_{\ell}(\hat
m)\right)\label{majo1}
\end{eqnarray}
that is
\begin{eqnarray} \|\tilde \ell-\ell\|_2^2
\leq \kappa_x^2 \|\ell-\ell_m^{(n)}\|_2^2
+2\kappa_x {\rm pen}_{\ell}(m) + x\kappa_x W_n(\hat m),\label{majo2}\end{eqnarray} where
\begin{equation}\label{wn} W_{n}(m')=
[\sup_{t\in B_{m, m'}(0,1)} |\nu_n(t)|^2-p_{\ell}(m, m')]_+.\end{equation}
The main point of the proof lies in studying $W_n(m')$, more precisely in finding
$p_{\ell}(m,m')$ such that
\begin{equation}\label{but}
{\mathbb E}(W_n(\hat m))  \leq \sum_{m'\in {\mathcal M}_{n,\ell}}
\mathbb{E}(W_n(m'))) \leq C/n,
\end{equation}  where $C$ is a constant.
In this case, combining \eref{majo2} and \eref{but} we infer that,
for all $m$ in ${\mathcal M}_{n, \ell}$,
$$
\mathbb{E}\|\ell-\tilde \ell\|_2^2 \leq \kappa_x^2 \|\ell-\ell_m^{(n)}\|_2^2 +
2\kappa_x{\rm pen}_{\ell}(m) + x\kappa_xC/n,$$  which can also be written
\begin{equation}\label{compromi}
\mathbb{E}\|\ell-\tilde \ell\|_2^2 \leq C_x\inf_{m\in {\mathcal
M}_{n, \ell} } \left[ \|\ell-\ell_m\|_2^2 + {\rm
pen}_{\ell}(m)\right] + C_xC'/n,
\end{equation}
where $C_x=\max(\kappa_x^2,2\kappa_x)$ suits, when $k_n\geq n$, and  \eref{mn}
and \eref{mn+} hold. It
remains thus to find $p_{\ell}(m,m')$ such that \eref{but} holds.

The process $W_n(m')$ is studied by using the
decomposition of $\nu_n(t) =\nu_{n,1}(t)+\nu_{n,2}(t)$ \mbox{ with }
\begin{eqnarray}
\label{nuni}
 \nu_{n,1}(t) = \frac 1n
\sum_{i=1}^n (f(X_i)u_t^*(Z_i) -\langle t, \ell\rangle)
\mbox{ and } \nu_{n,2}(t) = \frac 1n \sum_{i=1}^n \xi_iu_t^*(Z_i).
\end{eqnarray}
It follows that $W_n(m') \leq  2W_{n,1}(m')+ 2 W_{n,2}(m')$
where for $i=1,2$, \begin{equation}\label{wni} W_{n,i}(m')=
[\sup_{t\in B_{m, m'}(0,1)} |\nu_{n,i}(t)|^2-p_i(m, m')]_+,\mbox{
  and } p_{\ell}(m,m')=2p_1(m,m') + 2p_2(m,m').\end{equation}

\noindent $\bullet$ Study of $W_{n,1}$.\\
Since under \eref{fbornee}, $f$ is bounded on the support of $g$, we
apply a standard Talagrand's~(1996) inequality (see Lemma
\ref{lemtal} below that can be a fortiori applied to identically
distributed variables):
\begin{lem}\label{lemtal}
Let $U_1, \dots, U_n$ be independent random variables and
$\nu_n(r)=(1/n)\sum_{i=1}^n [r(U_i)-\mathbb{E}(r(U_i))]$
for $r$ belonging
to a countable class ${\mathcal R}$ of uniformly bounded
measurable functions. Then for $\epsilon>0$
\begin{equation}\label{talesp}
\mathbb{E}\left[\sup_{r\in {\mathcal
R}}|\nu_n(r)|^2-2(1+2\epsilon)H^2\right]_+ \leq \frac
6{K_1}\left(\frac vn e^{-K_1\epsilon \frac{nH^2}v}
 + \frac{8M_1^2}{K_1n^2C^2(\epsilon)}
e^{-\frac{K_1 C(\epsilon)\sqrt{\epsilon}}{\sqrt{2}}\frac{nH}{M_1}}\right),
\end{equation}
with $C(\epsilon)=\sqrt{1+\epsilon}-1$, $K_1$ is a universal constant, and where $$\sup_{r\in {\mathcal
R}}\|r\|_{\infty}\leq M_1, \;\;\;\; \mathbb{E}\left(\sup_{r\in
{\mathcal R}}|\nu_n(r)|\right)\leq H, \;\;\;\; \sup_{r\in
{\mathcal R}}\frac{1}{n}\sum_{i=1}^n{\rm Var}(r(U_i)) \leq v.$$
\end{lem}
The inequality (\ref{talesp}) is a straightforward consequence of
Talagrand's~(1996) inequality given in Ledoux~(1996) (or Birg\'e
and Massart~(1998)). Therefore
\begin{equation}\label{talstand}
\mathbb{E}[\sup_{t\in B_{m,m'}(0,1)}
|\nu_{n,1}(t)|^2-2(1+2\epsilon_1){\mathbb H}_1^2]_+ \leq
\kappa_1\left(\frac{v_1}n e^{-K_1 \epsilon_1 \frac{n{\mathbb
H}_1^2}{v_1}} + \frac{M_1^2}{n^2} e^{-K_2\sqrt{\epsilon_1}C(\epsilon_1) \frac{n{\mathbb
H_1}}{M_1}}\right),\end{equation}  where $K_2=K_1/\sqrt{2}$ and ${\mathbb H}_1$, $v_1$ and $M_1$
are defined by
${\mathbb E}(\sup_{t\in B_{m,m'}(0,1)}
|\nu_{n,1}(t)|^2) \leq {\mathbb H}_1^2$,
$$\sup_{t\in B_{m,m'}(0,1)} {\rm Var}(f(X_1)u_t^*(Z_1)) \leq v_1, \mbox{ and }
\sup_{t\in B_{m,m'}(0,1)} \|f(X_1)u_t^*(Z_1)\|_{\infty} \leq
M_1.$$ According to \eref{gamma} and \eref{var}, we propose to take \begin{equation}\label{m1}
M_1=M_1(m,m')=\kappa_{\infty, G}\sqrt{\lambda_1\Gamma(m^*)}.\end{equation}
For $v_1$, denoting by $P_{j,k}$, the quantity $P_{j,k}(m)= {\mathbb
E}\left[f^2(X_1)u_{\varphi_{m,j}}^*(Z_1)
u_{\varphi_{m,k}}^*(-Z_1)\right],$ write
\begin{eqnarray*}
\sup_{t\in B_{m,m'}(0,1)} {\rm Var}(f(X_1)u_t^*(Z_1)) &\leq &
( \sum_{j,k\in {\mathbb Z}} |P_{j,k}(m^*)|^2)^{1/2}.
\end{eqnarray*}
Arguing as in Comte\textit{ et al.}~(2005a),
let us define $\Delta_2(m,\Psi)$  by
\begin{eqnarray}\label{delta2} \Delta_2(m,\Psi)=D_m^2\int\int
\left|\frac{\varphi^*(x)\varphi^*(y)}{f_{\varepsilon}^*(D_m
x)f_{\varepsilon}^*(D_m y)} \Psi^*(D_m(x-y))\right|^2
dxdy\leq\lambda_2^2(\|\Psi\|_2) \Gamma_2^2(m^*),\nonumber\end{eqnarray}  with
\begin{equation}\label{gamma2}
 \Gamma_2(m^*)=  D_{m^*}^{2\alpha+
\min[(1/2-\rho/2),(1-\rho)]}\exp\{2\beta\sigma^\rho(\pi D_{m^*})^{\rho}\}
\end{equation} and $\lambda_2(\|\Psi\|_2)=\lambda_2(\alpha,\kappa_0,\beta,\sigma,\rho,
\|\Psi\|_2)$ given by
\begin{equation}\label{lambda2}
\lambda_2(\|\Psi\|_2)=\left\lbrace
\begin{array}{ll} \lambda_1(\alpha,\kappa_0,\beta,\sigma,\rho) &\mbox{ if }
\rho
>1,\\
\kappa_0^{-1}(2\pi)^{-1/2}
\lambda_1^{1/2}(\alpha,\kappa_0,\beta,\sigma,\rho)(1+\sigma^2\pi^2)^{\alpha/2}\|\Psi\|_2
&\mbox{ if } \rho\leq 1.\end{array}\right.\end{equation} Now, write $P_{j,k}$ as
\begin{eqnarray*}
P_{j,k}(m) &=& \int\!\!\!\!\int f^2(x)
u_{\varphi_{m,j}}^*(x+y)u_{\varphi_{m,k}}^*(-(x+y))
g(x)f_{\varepsilon}(y)dxdy \end{eqnarray*} that is
\begin{eqnarray*}
P_{j,k}(m) &=&\!\!\!D_{m}\!\int\!\!\!\!\int\!\!f^2(x)
\int\!\!\!\!\int\!e^{-i(x+y)uD_{m}} \frac{\varphi^*(u)
e^{iju}}{f_{\varepsilon}^*(D_{m}u)} e^{i(x+y)vD_{m}}
\frac{\varphi^*(v) e^{ikv}}{f_{\varepsilon}^*(D_{m}v)}dudv
g(x)f_{\varepsilon}(y) dxdy
\\
 &=&\!\!D_{m}\!\int\!\!\!\!\int
\frac{e^{iju+ikv}
\varphi^*(u)\varphi^*(v)}{f_{\varepsilon}^*(D_{m}u)f_{\varepsilon}^*(D_{m}v)}\left(
\int\!\!\!\!\int e^{-i(x+y)(u-v)D_{m}}
f^2(x)g(x)f_{\varepsilon}(y) dxdy\right) dudv \\ &=& D_{m}
\int\!\!\!\int \frac{e^{iju+ikv}
\varphi^*(u)\varphi^*(v)}{f_{\varepsilon}^*(D_{m}u)f_{\varepsilon}^*(D_{m}v)}
[(f^2g)*f_{\varepsilon}]^*((u-v)D_{m})dudv.
\end{eqnarray*}
By applying  Parseval's formula we get that $\sum_{j,k} |P_{j,k}(m)|^2$ equals
\begin{eqnarray*}
 D_{m}^2 \int\int \left|\frac{
\varphi^*(u)\varphi^*(v)}{f_{\varepsilon}^*(D_{m}u)f_{\varepsilon}^*(D_{m}v)}
[(f^2g)*f_{\varepsilon}]^*((u-v)D_{m})\right|^2 dudv =
\Delta_2(m, (f^2g)*f_{\varepsilon}).
\end{eqnarray*}
Since $\|(f^2g)*f_{\varepsilon}\|_2\leq \|f^2g\|_2\|f_{\varepsilon}\|_2=
{\mathbb E}^{1/2}(f^2(X_1))\|f_{\varepsilon}\|_2$, and
$\lambda_2(\|f^2g\|_2\|f_{\varepsilon}\|_2)\leq \mu_2$, by using the
definition of $\mu_2$ given in (\ref{mu1}),  we propose to take
\begin{equation}\label{v1} v_1=v_1(m,m')= \mu_2\Gamma_2(m^*).
\end{equation} Lastly, we have
${\mathbb E}[\sup_{t\in B_{m,m'}(0,1)}
|\nu_{n,1}(t)|^2]
\leq {\mathbb
E}(f^2(X_1))\lambda_1\Gamma(m^*)/n$  and thus we propose to
take \begin{eqnarray}
 \label{p1} {\mathbb
H}_1^2={\mathbb
H}_1^2(m,m')= {\mathbb
E}(f^2(X_1)) \lambda_1\Gamma(m^*)/n.\end{eqnarray} It follows from (\ref{talstand}),
(\ref{m1}), (\ref{v1}) and (\ref{p1}) that if
$$p_1(m,m')= 2(1+2\epsilon_1){\mathbb H}_1^2=2(1+2\epsilon_1){\mathbb
E}(f^2(X_1)) \lambda_1\Gamma(m^*)/n$$ then
\begin{eqnarray}\mathbb{E}(W_{n,1}(m'))&\leq&\mathbb{E}\left[\sup_{t\in B_{m,m'}(0,1)}
|\nu_{n,1}(t)|^2-2(1+2\epsilon_1){\mathbb H}_1^2\right]_+ \leq A_1(m^*)+
B_1(m^*)
\end{eqnarray}
with
\begin{eqnarray}
\label{a1}
A_1(m)&=&K_3\frac{\mu_2\Gamma_2(m)}n
\exp\left(-K_1\epsilon_1{\mathbb
E}(f^2(X_1))\frac{\lambda_1\Gamma(m)}{\mu_2\Gamma_2(m)}\right)\\\label{b1}
\mbox{ and }
B_1(m)&=&K_3\frac{\kappa_{\infty, G}^2\lambda_1\Gamma(m)}{n^2}\exp\left\lbrace-K_2\sqrt{\epsilon_1}C(\epsilon_1)\frac{\sqrt{{\mathbb E}(f^2(X_1))}}{\kappa_{\infty, G}}\sqrt{n}\right\rbrace.
\end{eqnarray}
Since $\forall m\in {\mathcal M}_{n, \ell}$, $\Gamma(m)\leq n$ and
$|{\mathcal M}_{n,\ell}|\leq n$, there exist some
constants $K_4$ and $c$ such that
$$\sum_{m\in {\mathcal M}_{n, \ell}} B_1(m^*) \leq K_3 \|f\|_{\infty,
G}^2\lambda_1\exp[-K_4\sqrt{{\mathbb
E}(f^2(X_1))}\sqrt{n}/\kappa_{\infty, G}]\leq c/n.$$
Let us now come to the study of $A_1(m^*)$.

\paragraph{\textbf{1) Case $0\leq \rho <1/3$}} In that case,
$\rho\leq (1/2-\rho/2)_+$ and the choice $\epsilon_1=1/2$
ensures the convergence of $\sum_{m'\in {\mathcal M}_{n, \ell}}
A_1(m^*)$. Indeed, if we denote by $\psi=2\alpha+
\min[(1/2-\rho/2),(1-\rho)]$, $\omega=(1/2-\rho/2)_+$,
$K'=\kappa_2\lambda_1/\mu_2$, then for $a,b\geq 1$,  we infer that
\begin{eqnarray}\nonumber
\max(a, b)^{\psi}e^{2\beta\sigma^\rho\pi^{\rho} \max(a,b)^{\rho}}e^{-K'\xi^2\max(a,b)^{\omega}}\!\!\!\!\!\!&\leq &
(a^{\psi}e^{2\beta\sigma^\rho\pi^{\rho} a^{\rho}}+b^{\psi}
e^{2\beta\sigma^\rho\pi^{\rho}
b^{\rho}})e^{-(K'\xi^2/2)(a^{\omega} +
b^{\omega})}
\end{eqnarray}
is bounded by
\begin{eqnarray}\label{eqmax}
a^{\psi}e^{2\beta\sigma^\rho\pi^{\rho}
a^{\rho}}e^{-(K'\xi^2/2)a^{\omega}} e^{-(K'\xi^2/2)
b^{\omega}}+b^{\psi} e^{2\beta\sigma^\rho\pi^{\rho} b^{\rho}}
e^{-(K'\xi^2/2) b^{\omega})}.\end{eqnarray} Since the function
$a\mapsto a^{\psi}e^{2\beta\sigma^\rho\pi^{\rho}
a^{\rho}}e^{-(K'\xi^2/2)a^{\omega}} $ is bounded on ${\mathbb
R}^+$ by a constant, depending on $\alpha$, $\rho$ and $K'$
only, and since $Ak^{\rho}-\beta k^{\omega} \leq -(\beta/2)k^{\omega}$
for any $k\geq 1$, it follows that $\sum_{m'\in {\mathcal M}_{n,
\ell}}A_1(m^*)\leq C/n.$

\paragraph{\textbf{2) Case $\rho=1/3$}.}
In that case, $\rho =
(1/2-\rho/2)_+$, and $\omega=\rho$. We
choose $\epsilon_1=\epsilon_1(m,m')$ such that $2\beta\sigma^\rho\pi^{\rho} D_{m^*}^{\rho}
- K'\mathbb{E}(f^2(X_1))    \epsilon_1 D_{m^*}^{\rho}= -2\beta\sigma^\rho\pi^{\rho} D_{m^*}^{\rho}$ that
is, since $K'=K_1\lambda_1/\mu_2$,
$\epsilon_1=\epsilon_1(m,m')=(4\beta\sigma^\rho\pi^{\rho}\mu_2)/(K_1\lambda_1\mathbb{E}(f^2(X_1))).$
\paragraph{\textbf{3) Case $\rho>1/3$}.}
In that case, $\rho >
(1/2-\rho/2)_+$. Bearing in mind the inequality (\ref{eqmax}) we
choose $\epsilon_1=\epsilon_1(m,m')$ such that $2\beta\sigma^\rho\pi^{\rho} D_{m^*}^{\rho}
- K'\mathbb{E}(f^2(X_1))    \epsilon_1 D_{m^*}^{\omega}= -2\beta\sigma^\rho\pi^{\rho} D_{m^*}^{\rho}$ that
is, since $K'=K_1\lambda_1/\mu_2$,
$\epsilon_1=\epsilon_1(m,m')=(4\beta\sigma^\rho\pi^{\rho}\mu_2)/(K_1\lambda_1\mathbb{E}(f^2(X_1)))D_{m^*}^{\rho-\omega}.$

These choices ensure that $\sum_{m'\in {\mathcal M}_{n,
\ell}}A_1(m^*)$ is less than $C/n$.

\noindent $\bullet$ Study of $W_{n,2}$.\\
Denote by
\begin{eqnarray}
\label{Hxi}
{\mathbb H}^2_{\xi}(m,m') = (n^{-1}\sum_{i=1}^n
\xi_i^2) \lambda_1\Gamma(m^*)/n,\end{eqnarray} with $( n^{-1}\sum_{i=1}^n \xi_i^2)\lambda_1\Gamma(m)/n
=  (n^{-1} \sum_{i=1}^n \xi_i^2-\sigma_{\xi}^2)
\lambda_1\Gamma(m)/n + \sigma^2_{\xi} \lambda_1\Gamma(m)/n$ bounded by
\begin{eqnarray*}
(n^{-1}\sum_{i=1}^n
\xi_i^2-\sigma_{\xi}^2)\1_{\{n^{-1} |\sum_{i=1}^n
(\xi_i^2-\sigma^2_{\xi})|
>\sigma_{\xi}^2/2\}} \lambda_1\Gamma(m)/n  +
3\sigma^2_{\xi}\lambda_1\Gamma(m)/(2n).\end{eqnarray*}
Consequently
${\mathbb H}_{\xi}^2(m,m') \leq {\mathbb H}_{\xi,1}(m,m')+ {\mathbb H}_{\xi,2}(m,m')$
where  $${\mathbb H}_{\xi,1}(m,m') = ( n^{-1}\sum_{i=1}^n
\xi_i^2-\sigma_{\xi}^2)\1_{\{n^{-1}| \sum_{i=1}^n
\xi_i^2-\sigma^2_{\xi}|
>\sigma_{\xi}^2/2\}} \lambda_1\Gamma(m^*)/n \mbox{ and }{\mathbb H}_{\xi,2}(m,m')=
3\sigma^2_{\xi}\lambda_1\Gamma(m^*)/(2n).$$
By applying \eref{nuni}
we infer that ${\mathbb E}[\sup_{t\in
B_{m, m'}(0,1)} |\nu_{n,2}(t)|^2 -p_2(m, m')]_+$ is bounded by
\begin{multline*} {\mathbb E}[2\!\!\!\sup_{t\in B_{m, m'}(0,1)}
(n^{-1} \sum_{i=1}^n \xi_i (u_t^*(Z_i) - \langle t,
g)\rangle)^2- 4(1+2\epsilon_2){\mathbb H}_{\xi}^2(m,m')]_+
+ 2\|g\|^2_2 {\mathbb E}[(n^{-1}\sum_{i=1}^n \xi_i)^2] \\+ {\mathbb
E}[4(1+2\epsilon_2){\mathbb H}_{\xi}^2(m,m') - p_2(m,m')]_+,
\end{multline*}
that is
\begin{multline}\label{decompprel}
{\mathbb E}[\sup_{t\in
B_{m, m'}(0,1)} |\nu_{n,2}(t)|^2 -p_2(m, m')]_+
\\
\leq 2{\mathbb E}[\sup_{t\in B_{m, m'}(0,1)} (n^{-1}\sum_{i=1}^n \xi_i (u_t^*(Z_i) - \langle t, g\rangle ) )^2-
2(1+2\epsilon_2){\mathbb H}_{\xi}^2(m,m')]_++
2{\|g\|^2_2\sigma_{\xi}^2}/n\\
+ 4(1+2\epsilon_2){\mathbb E}|{\mathbb H}_{\xi,1}(m,m')| +
{\mathbb E}[4(1+2\epsilon_2) {\mathbb H}_{\xi,2}(m,m') - p_2(m,m')]_+.
\end{multline}
Since we only consider dimensions $D_m$ such that $\Gamma(m)/n$ is bounded by
some constant $\kappa$, we get that for some $p\geq 2$, ${\mathbb E}|{\mathbb
  H}_{\xi,1}(m,m')|$ is bounded by
\begin{eqnarray*}
\kappa \lambda_1{\mathbb E}[|\frac{1}n\sum_{i=1}^n
\xi_i^2-\sigma_{\xi}^2|\1_{\{n^{-1} |\sum_{i=1}^n
(\xi_i^2-\sigma^2_{\xi})|>\sigma_{\xi}^2/2\}}]
\leq {\kappa \lambda_1 2^{p-1}}{\mathbb E}[|n^{-1}\sum_{i=1}^n
\xi_i^2-\sigma_{\xi}^2|^p]/{\sigma_{\xi}^{2(p-1)}}
\end{eqnarray*}
According to Rosenthal's inequality (see Rosenthal~(1970)), we find that, for $\sigma^p_{\xi,p}:= {\mathbb E}(|\xi|^p), \;\;
\sigma_{\xi,2}^2=\sigma_{\xi}^2$, $${\mathbb E}|n^{-1}\sum_{i=1}^n
\xi_i^2-\sigma_{\xi}^2|^p\leq C'(p)\left(
\sigma_{\xi,2p}^{2p}n^{1-p} +
\sigma_{\xi,4}^{2p}n^{-p/2}\right).$$
 Now, the assumption \eref{TFfeps}
implies that $\alpha>1/2$, therefore $|\mathcal{M}_n|\leq
\sqrt{n}$ and consequently, by choosing $p=3$ this leads to
$\sum_{m' \in \mathcal{M}_n}{\mathbb E}|{\mathbb H}_{\xi,1}(m,m')|\leq C(\sigma_{\xi,6},\sigma_{\xi})/n.$
The last term of the inequality (\ref{decompprel}) vanishes
as soon as
\begin{equation}\label{p2m}
p_2(m,m')= 4(1+2\epsilon_2) {\mathbb H}_{\xi,2}(m,m')
=6(1+2\epsilon_2)\lambda_1\sigma^2_{\xi}\Gamma(m^*)/n.\end{equation}
For this choice of $p_2(m,m')$,
the inequality (\ref{decompprel}) becomes ${\mathbb E}[\sup_{t\in B_{m,
\hat m}(0,1)} |\nu_{n,2}(t)|^2-p_2(m, \hat m)]_+$ is less than
\begin{multline*}
2\sum_{m'\in {\mathcal M}_{n, \ell}} {\mathbb E}[\sup_{t\in B_{m,m'}(0,1)}
(n^{-1} \sum_{i=1}^n \xi_i (u_t^*(Z_i) - \langle t,
g\rangle ))^2- 2(1+2\epsilon_2){\mathbb H}_{\xi}^2(m,m')]_+\\
 +  2\|g\|^2_2\sigma_{\xi}^2/n + 4C(1+2\epsilon_2)/n.
\end{multline*}
Then we apply the following Lemma to reach the same kind of result
as \eref{talstand} for $W_{n,1}$.
\begin{lem}\label{talcond}
Under the assumptions of Theorem \ref{adaptatifgl}, if $\mathbb{E}\vert \xi_1\vert^6<\infty$, then for some given $\epsilon_2>0$:
\begin{multline} \label{inegnew}   \sum_{m'\in {\mathcal M}_{n,\ell}} {\mathbb
E}\left[\sup_{t\in B_{m, m'}(0,1)} \left(\frac 1n \sum_{i=1}^n
\xi_i (u_t^*(Z_i) - \langle t,
g\rangle \right)^2- 2(1+2\epsilon_2){\mathbb H}_{\xi}^2(m,m')\right]_+ \\
\leq K_1\left\{\sum_{m'\in {\mathcal M}_{n,\ell}} \left[\frac{\sigma_{\xi}^2\mu_2\Gamma_2(m^*)}n
\exp\left( -K_1\epsilon_2
\frac{\lambda_1\Gamma(m^*)}{\mu_2\Gamma_2(m^*)}\right)\right] +
\left(1+ \frac{\ln^4(n)}{\sqrt{n}}\right) \frac 1{n} \right\},
\end{multline}
where $\mu_2$ and $\Gamma_2(m)$ are defined by (\ref{mu1}) and
(\ref{gamma2}) and $K_1$ is a constant depending on the moments of
$\xi$. The constant $\mu_2$ can be replaced by $\lambda_2(\|h\|_2)$ where $\lambda_2$ is defined
by (\ref{lambda2}).
\end{lem}
\noindent By analogy with (\ref{a1}) we denote by
\begin{equation}\label{a2} A_2(m^*)=
\frac{K_1\sigma_{\xi}^2}n \mu_2\Gamma_2(m^*) \exp\left( -K_1\epsilon_2
\frac{\lambda_1\Gamma(m^*)}{\mu_2\Gamma_2(m^*)}\right)=
\frac{K_1\sigma_{\xi}^2\mu_2\Gamma_2(m^*)}{n} \exp\left(-\kappa_2\epsilon_2 \frac{\lambda_1}{\mu_2}D_{m^*}^{(1/2-\rho/2)_+}\right).
\end{equation}
With  $p_2(m,m')$ given
by (\ref{p2m}), by gathering (\ref{talstand}) and  (\ref{inegnew}), we find, for $W_{n,2}$ defined by (\ref{wni}),
$$\mathbb{E}(W_{n,2}(\hat m))\leq K\sum_{m' \in
\mathcal{M}_n}A_2(m^*) +C(
1+\ln(n)^6/n)/n + K'/n .$$ The sum $\sum_{m' \in
\mathcal{M}_n}A_2(m^*)$ is bounded in the same way as the sum
$\sum_{m' \in \mathcal{M}_n}A_1(m^*)$ with
$\epsilon_2=\epsilon_1=1/2$ if $0\leq \rho< 1/3$ and
$\epsilon_1(m,m')$ replaced by
$\epsilon_2=\epsilon_2(m,m')=\mathbb{E}(f^2(X_1))\epsilon_1(m,m'),$
when $\rho\geq  1/3$ that is
$\epsilon_2(m,m')=(4\beta\sigma^\rho\pi^\rho\mu_2)/(K_1\lambda_1)D_{m^*}^{\rho-\omega}$.
These choices ensure that $\sum_{m'\in {\mathcal M}_{n,
\ell}}A_2(m^*)$ is less than $C/n$. The result follows by taking as
announced in \eref{wni}, $p_\ell(m,m')= 2p_1(m,m')+ 2p_2(m,m')$, that is
$p_\ell(m,m')=
4[(1+2\epsilon_1(m,m')){\mathbb E}(f^2(X_1))
+ 3(1 +2\epsilon_2(m,m'))\sigma_{\xi}^2]
\lambda_1\Gamma(m^*)/n,$
and more precisely if $0\leq\rho < 1/3$,
\begin{eqnarray}
 p_\ell(m,m')=
24 \mathbb{E}(Y_1^2)\lambda_1\Gamma(m^*)/n,
\label{pl1}
\end{eqnarray}
and if $\rho\geq 1/3$,\begin{eqnarray}
 p_\ell(m,m')
=4[3\mathbb{E}(Y_1^2)+{32\beta\sigma^\rho\pi^\rho\mu_2}D_{m^*}^{\rho-\omega}/{k_1\lambda_1}]\lambda_1\Gamma(m^*)/n.
\label{pl3}
\end{eqnarray}
Consequently if $0\leq \rho<1/3$, we take
$\mbox{pen}_{\ell}(m)=\kappa{\mathbb E}(Y_1^2)
\lambda_1\Gamma(m)/n
,$
and if $\rho \geq 1/3$ we take
$
\mbox{pen}_{\ell}(m)=\kappa[\mathbb{E}(Y_1^2)+\beta\sigma^\rho\pi^\rho\mu_2D_m^{\rho-\omega}/{k_1\lambda_1}]\lambda_1\Gamma(m)/n,
$ for some numerical constants $\kappa$. Note that for $\rho=1/3$, $\rho-\omega=0$ and the second penalty has the same order as the first one with
a different multiplicative constant. \hfill $\Box$

\subsection{Proof of Lemma \protect \ref{talcond}, by using a
conditioning argument}

We work conditionally to the $\xi_i$'s and ${\mathbb E}_{\xi}$ and
${\mathbb P}_{\xi}$ denote the conditional expectations and
probability for fixed $\xi_1, \dots, \xi_n$.

We apply Lemma \ref{lemtal} with $f_t(\xi_i,Z_i)=\xi_iu_t^*(Z_i)$, conditionally to the $\xi_i$'s to the
random variables $(\xi_1,Z_1), \dots,$ $
(\xi_n, Z_n)$ which are
independent but non identically distributed since the $\xi_i$'s
are fixed constants.
Let $Q_{j,k}={\mathbb
E}[u_{\varphi_{m,j}}^*(Z_1)u_{\varphi_{m,k}}^*(-Z_1)].$
Straightforward calculations give that for $\mathbb{H}_\xi(m,m')$ defined in
\eref{Hxi} we have
$${\mathbb E}_{\xi}^2[\sup_{t\in B_{m,m'}(0,1)}n^{-1}
\sum_{i=1}^n \xi_i(u_t^*(Z_i)-\langle t,g\rangle )] \leq
\mathbb{H}_\xi^2(m,m').
$$
Again, arguing as in
Comte\textit{ et al.}~(2005a), $\sum_{j,k}|Q_{j,k}|^2\leq
\Delta_2(m,h)\leq \lambda_2(\|h\|_2)\Gamma_2(m, \|f_{\varepsilon}\|_2)$ with
$\|h\|_2\leq \|f_{\varepsilon}\|_2$, where
$\Delta_2(m,h)$ is defined by (\ref{delta2}), $\lambda_2$ by
(\ref{lambda2}), $\Gamma_2(m)$ by (\ref{gamma2}), $\mu_2$ by
(\ref{mu1}). We now write that\begin{eqnarray*} \sup_{t\in
B_{m,m'}(0,1)}(n^{-1} \sum_{i=1}^n {\rm
Var}_{\xi}(\xi_iu_t^*(Z_i)))
&\leq &
(n^{-1}\sum_{i=1}^n \xi_i^2)\mu_2\Gamma_2(m^*,
\|f_{\varepsilon}\|_2)\end{eqnarray*}
and thus we take
\begin{eqnarray*}v_{\xi}(m,m')=(n^{-1}\sum_{i=1}^n \xi_i^2)\mu_2\Gamma_2(m^*,
\|f_{\varepsilon}\|_2).\end{eqnarray*}   Lastly, since
$$\sup_{t \in B_{m,m'}(0,1)} \|f_t\|_{\infty} \leq 2 \max_{1\leq i\leq
n}|\xi_i| \sqrt{\Delta(m^*)}\leq 2\max_{1\leq i\leq n}|\xi_i|
\sqrt{\lambda_1\Gamma(m^*)}$$
we take $M_{1,\xi}(m,m')=2\max_{1\leq i\leq n}|\xi_i|
\sqrt{\lambda_1\Gamma(m^*)}.$ By applying Lemma \ref{lemtal}, we get
for some constants $\kappa_1$, $\kappa_2$, $\kappa_3$
\begin{multline*}
{\mathbb E}_{\xi}[\sup_{t\in B_{m,m'}(0,1)}\nu_{n,1}^2(t) -
2(1+2\epsilon){\mathbb H}_{\xi}^2]_+  \leq
\displaystyle\kappa_1 \left[ \frac{\mu_2\Gamma_2(m^*)}{n^2}(\sum_{i=1}^n
\xi_i^2) \exp\left\lbrace{-\kappa_2
\epsilon\frac{\lambda_1\Gamma(m^*)}{\mu_2\Gamma_2(m^*)}}\right\rbrace\right.\\ \left.+
\frac{\lambda_1\Gamma(m^*)}{n^2}(\max_{1\leq i\leq n}\xi_i^2)
 \exp\left\lbrace{-\kappa_3\sqrt{\epsilon}C(\epsilon) \frac{\sqrt{\sum_{i=1}^n
\xi_i^2}}{\max_i|\xi_i|}}\right\rbrace \right]\end{multline*} To relax the
conditioning, it suffices to integrate with respect to the law of
the $\xi_i$'s the above expression. The first term in the bound
simply becomes:
$$\sigma_{\xi}^2 \mu_2\Gamma_2(m^*)
\exp[-\kappa_2 \epsilon\lambda_1\Gamma(m^*)/(\mu_2
\Gamma_2(m^*)])/n$$ and has the same order as in the case of
bounded variables. The second term is bounded by
\begin{equation}\label{secondterme} \frac{\lambda_1 \Gamma(m^*)}{n^2} {\mathbb
E}\left[(\max|\xi_i|^2)\exp\left(-\kappa_3\sqrt{\epsilon}C(\epsilon)
\frac{\sqrt{\sum_{i=1}^n \xi_i^2}}{\max_{1\leq i\leq n}
|\xi_i|}\right)\right].\end{equation}

Since we only consider dimensions $D_m$ such that the penalty term
is bounded, we have $\Gamma(m)/n\leq K$ and the sum of the above
terms for $m\in {\mathcal M}_{n, \ell}$ and $|{\mathcal M}_{n,
\ell}|\leq n$ is less than
$$ \lambda_1 {\mathbb
E}\left[\left(\max_{1\leq i\leq n} \xi_i^2\right)
\exp\left(-\kappa_3\sqrt{\epsilon}C(\epsilon) \frac{\sqrt{\sum_{i=1}^n
\xi_i^2}}{\max_{1\leq i\leq n} |\xi_i|}\right)\right].$$ We need
to study when such a term is less than $c/n$ for some constant
$c$. We bound $\max_i|\xi_i|$ by $b$ on the set $\{\max_i |\xi_i| \leq
b\}$ and the exponential by 1 on the set $\{\max_i |\xi_i| \geq b\}$ and
by denoting $\mu_\epsilon=\kappa_3\sqrt{\epsilon}C(\epsilon)$, this yields
\begin{eqnarray*} && {\mathbb E}\left[\max_{1\leq i\leq n} \xi_i^2
\exp\left(-\mu_\epsilon \sqrt{\frac{\sum_{i=1}^n \xi_i^2}{\max_{1\leq i\leq
n}\xi_i^2}}\right)\right] \\
&\leq & b^2 {\mathbb E}\left(\exp(-\mu_\epsilon
\frac{\sqrt{\sum_{i=1}^n \xi_i^2}}b) \right) +{\mathbb
E}\left(\max_{1\leq i\leq n} \xi_i^2\1_{\{\max_{1\leq i\leq
n}|\xi_i|\geq b\}}\right)\\
&\leq &  b^2\left[{\mathbb
E}\left(\exp(-\mu_\epsilon \sqrt{n\sigma_{\xi}^2/(2b^2)}\right) + {\mathbb
P}\left(|\frac 1n\sum_{i=1}^n \xi_i^2-\sigma_{\xi}^2|\geq
\sigma^2_{\xi}/2\right)\right] +b^{-r}{\mathbb E}(\max_{1\leq
i\leq n}|\xi_i|^{r+2})\\
&\leq & b^2e^{-\mu_\epsilon
\sqrt{n}\sigma_{\xi} / (\sqrt{2}b)} + b^2 2^p\sigma_{\xi}^{-2p} {\mathbb E}\left(\left|\frac
1n\sum_{i=1}^n \xi_i^2-\sigma_{\xi}^2\right|^p\right)+
b^{-r}{\mathbb E}(\max_{1\leq i\leq n}|\xi_i|^{r+2}).
\end{eqnarray*}
Again by applying Rosenthal's inequality (see Rosenthal~(1970)), we get that
\begin{multline*}{\mathbb E}\left[\max_{1\leq i\leq n} \xi_i^2
\exp\left(-\mu_\epsilon \sqrt{\frac{\sum_{i=1}^n \xi_i^2}{\max_{1\leq i\leq
n}\xi_i^2}}\right)\right]\\
\leq b^2 e^{-\mu_\epsilon \sqrt{n}\sigma_{\xi} / (\sqrt{2}b)} + b^2\frac{2^p}{\sigma_{\xi}^{2p}} \frac{C(p)}{n^p}  [ n{\mathbb
E}(|\xi_1^2-\sigma_{\xi}^2|^p) + (n{\mathbb
E}(\xi_1^4))^{p/2}] + n {\mathbb E}(|\xi_1|^{r+2})b^{-r}
\end{multline*} also bounded by
\begin{eqnarray*} b^2e^{-\mu_\epsilon \sqrt{n}\sigma_{\xi} / (\sqrt{2}b)} + C'(p)b^2
\sigma_{\xi, 2p}^{2p} 2^p\sigma_{\xi}^{-2p}
[ n^{1-p} + n^{-p/2}] + n\sigma_{\xi,
r+2}^{r+2}b^{-r}.
\end{eqnarray*}
Since $\mathbb{E}\vert \xi_1\vert^6<\infty$, we take
$p=3$, $r=4$, $b=\sigma_{\xi}\sqrt{\epsilon} C(\epsilon)\kappa_3\sqrt{n}/[2\sqrt{2}
(\ln(n)-\ln\ln n)]$ and for any $n \geq 3$, and for $C_1$ and
$C_2$ some constants
depending on the moments of $\xi$, we find that
$${\mathbb E}\left\{\left(\max_{1\leq i\leq n}\xi_i^2\right)
\exp\left(-\kappa_3 \sqrt{\epsilon} C(\epsilon) \sqrt{\sum_{i=1}^n
\xi_i^2/\max_{1\leq i\leq n}\xi_i^2}\right)\right\} \leq
\frac{C_1}{\sqrt{n}} +  C_2 \left(\frac{\ln^4(n)}{\sqrt{n}}
\right) \frac 1{\sqrt{n}}.$$ Then the sum over ${\mathcal
M}_{n, \ell}$ with cardinality less than $\sqrt{n}$ of the terms
in (\ref{secondterme}) is bounded by $C(1+\ln(n)^4/\sqrt{n})/n$
for some constant $C$, by using again that $\Gamma(m^*)/n$ is
bounded.

\subsection{Proof of Theorem \ref{adaptf}}

Let $\tilde E_n$ be the event $\tilde E_n=\{\parallel g-\tilde g \parallel_{\infty,A} \leq
g_0/2\}$.
Since $g(x)\geq g_0$ for any $x$ in $A$, then, on $\tilde E_n$,  $\tilde g (x)\geq g_0/2$ also
for any $x$ in $A$. It follows that
\begin{eqnarray}
\label{eq21}\mathbb{E}\|(f-\tilde f)\1_A\1_{\tilde E_n}\|_2^2\leq
{8}{g_0^{-2}}\mathbb{E}\|\tilde \ell-\ell\|_2^2
+{8\|\ell\|_{\infty,A}^2}{g_0^{-4}}\mathbb{E}\|\tilde g-g\|_2^2,
\end{eqnarray}
where $\|\ell\|_{\infty,A}\leq g_1\kappa_{\infty,G}$.
Using that $\|\tilde f\|_{\infty,A}\leq a_n$, we obtain
\begin{eqnarray}
\label{eq22}\mathbb{E}[\|(f-\tilde f)\1_A\1_{\tilde E_n^c}\|_2^2]\leq 2(a_n^2+\|f\|_{\infty,A}^2)\lambda(A)\mathbb{P}(\tilde E_n ^c),\end{eqnarray} where $\lambda(A)=\int _A dx$.
It follows that for $\hat m_\ell=\hat m_\ell(n)$, $\hat m_g=\hat m_g(n)$, if
$a_n\mathbb{P}(\tilde E_n^c)=o(n^{-1})$, then \eref{doubleinf} is proved by applying Theorem \ref{adaptatifgl}.
We now come to the study of $\mathbb{P}(\tilde E_n^c)$ by writing that
$
\mathbb{P}(\tilde E_n^c)=\mathbb{P}\left(\|g-\tilde g\|_{\infty}>g_0/2\right)
=\mathbb{P}\left(\|g-g^{(n)}_{\hat m_g}+g^{(n)}_{\hat m_g}-\tilde g\|_{\infty}>g_0/2\right).$
By applying Lemma \ref{connectnorm}:
\begin{lem}
\label{connectnorm}
Let $g$ belongs to $\mathcal{
  S}_{a_g,\nu_g,B_g}(\kappa_{a_g})$ defined by \eref{super} with
$a_g>1/2$. Then for $t\in S_m$, $ \|t\|_{\infty}\leq \sqrt{D_m}\|t\|_2$ and
$\|g-g_m\|_{\infty} \leq (2\pi)^{-1}\sqrt{\pi D_m} ((\pi D_m)^2+1)^{-a_g/2}
\exp(-B_g |\pi D_m|^{r_g}) A_g^{1/2}$.
\end{lem}
\noindent and by arguing as for  $\parallel
\ell_m-\ell_m^{(n)}\parallel_2^2$, we get that $\|g-g^{(n)}_{\hat m_g}\|_\infty\leq
\|g-g_{\hat m_g}\|_\infty+\|g_{\hat m_g}-g^{(n)}_{\hat m_g}\|_\infty$ also
bounded by
\begin{eqnarray*}
\sqrt{\kappa (\kappa_{\mathcal{G}}+1)} {D_{\hat m_g}^{3/2}}/\sqrt{k_n}}
+(2\pi)^{-1}{\sqrt{\pi D_{\hat m_g}} ((\pi D_{\hat m_g})^2+1)^{-a_g/2}
\exp(-B_g |\pi D_{\hat m_g}|^{\nu_g}) A_g^{1/2}.\end{eqnarray*}
Consequently, $\|g-g^{(n)}_{\hat m_g}\|_\infty$
tends to zero as soon as $g$ belongs to some space $\mathcal{ S}_{a_g,\nu_g,B_g}(\kappa_{a_g})$ defined by
\eref{super} with $a_g>1/2$ if $r_g=0$ and since $k_n\geq n^{3/2}$ and $D_{\hat m_g}=o(\sqrt{n})$ for 
$\alpha>1/2$. It follows that for $n$ large enough, $\|g-g^{(n)}_{\hat m_g}\|_\infty\leq g_0/4$ and consequently
$\mathbb{P}(\tilde E_n^c)\leq
\mathbb{P}[\|g^{(n)}_{\hat m_g}-\tilde g\|_{\infty}>g_0/4].
$
By applying again Lemma \ref{connectnorm}, since $g^{(n)}_{\hat m_g}-\tilde g$ belongs to
$S_{\hat m_g}$, we get that
\begin{eqnarray}
\mathbb{P}(\tilde E_n^c)\leq
\mathbb{P}[\|g^{(n)}_{\hat m_g}-\tilde g\|_2>g_0/(4\sqrt{D_{\hat m_g}})].
\end{eqnarray}
In this context, we have
\begin{eqnarray}
\label{bsupchap}
\|g^{(n)}_{\hat m_g}-\tilde g_{\hat m_g}\|_2^2&=&\sum_{\vert j\vert \leq k_n}(\hat
a_{\hat m_g,j}-a_{\hat m_g,j})^2=
\sum_{\vert j\vert \leq k_n}\nu_{n,g}^2(\varphi_{\hat m_g,j})=\sup_{t \in B_{\hat m_g}(0,1)}\nu_{n,g}^2(t).
\end{eqnarray}
Consequently,
\begin{eqnarray*}
\mathbb{P}(\tilde E_n^c)&\leq& \mathbb{P}[\sup_{t \in B_{\hat m_g}(0,1)}|\nu_{n,g}(t)|\geq
  g_0/(4\sqrt{D_{\hat m_g}})]
\leq  \sup_{m\in \mathcal{M}_n}\mathbb{P}[\sup_{t \in B_{\hat m_g}(0,1)}|\nu_{n,g}(t)|\geq
  g_0/(4\sqrt{D_{m}})]\\
&\leq& \sum_{m\in \mathcal{M}_n}\mathbb{P}[\sup_{t \in B_{\hat m_g}(0,1)}|\nu_{n,g}(t)|\geq
  g_0/(4\sqrt{D_{m}})].
\end{eqnarray*}
We apply Talagrand's~(1996) inequality as given
in Birg\'e and Massart~(1998), to
get that if we take $\lambda= g_0/(8\sqrt{D_{m}}) $ and if we
ensure $2H<g_0/(8\sqrt{D_{m}})$, then $\mathbb{P}[ \sup_{t\in
B_{m}(0,1)}|\nu_{n,g}(t)|\geq g_0/(4\sqrt{D_{m}})]\leq
3\exp\left[-K_1'n\left(\min[(D_{m}v)^{-1}, (M_1\sqrt{D_{m}})^{-1}]
\right) \right].$ This yields
\begin{eqnarray}
\mathbb{P}(\tilde E_n^c)&\leq& K\sum_{m\in \mathcal{M}_n}\{
\exp[-K_1'n/(M_1\sqrt{D_{m}})]+\exp[-K_1'n/(D_{m}v)]\}.
\end{eqnarray}
Since we only consider $D_m$ such that $D_{m}\leq \sqrt{n}$,
\begin{eqnarray*}
a_n|\mathcal{M}_n|\exp[-K_1'n/(M_1\sqrt{D_{m}})]\leq
a_n|\mathcal{M}_n|\exp(-K"n^{1/4})=o(n^{-1}).
\end{eqnarray*}
We only consider $D_m$ such that
$\Gamma(m)/n$ tends to zero.
Consequently, when $\rho>0$ then $D_m\leq \left(\ln n/(2\beta\sigma^\rho+1)\right)^{1/\rho}$ which combined with the fact that $v\leq
D_{m}^{2\alpha+1-\rho}\exp(2\beta\sigma^\rho\pi^\rho  D_{m}^\rho)$
gives that $a_n|\mathcal{M}_n|\exp\left(-K_1'n/(D_{m}v)\right)=o(1/n).$

When $\rho=0$, then $v=\mu_1 D_{m}^{2\alpha+1/2}$ and consequently, as $D_m\leq (n/\ln(n))^{1/(2\alpha +1)}\leq n^{1/(2\alpha+1)}$,
$$\exp(-K_1'n/(D_{m}v))\leq
\exp(-K"n/(D_m^{2\alpha+3/2}))\leq \exp( -K" n^{1/(4(\alpha +1))}).$$ Analogously, $\sqrt{D_m}H\leq 1/\sqrt{\ln(n)}$ in the worst case corresponding to $\rho=0$, for $D_m\leq (n/\ln(n))^{1/(2\alpha+2)}$, tends to zero and therefore is bounded by $g_0/8$ for $n$ great enough.
We conclude that if we only consider $D_m$ such that $D_m\leq n^{1/(2\alpha+2)}$ then $a_n\mathbb{P}(\tilde E_n^c)=o(1/n),$ and the result follows by applying the inequalities \eref{eq21} and \eref{eq22}.
\hfill $\Box$\\

\noindent {\bf Proof of Lemma \ref{connectnorm}.}
For $t \in S_m$, written as $t(x)=\sum_{j\in {\mathbb
Z}} \langle t, \varphi_{m,j}\rangle \varphi_{m,j}(x)$ and $|t(x)|^2 \leq \sum_{j\in{\mathbb Z}}|\langle t, \varphi_{m,j}\rangle |^2\sum_{j\in{\mathbb
Z}}|(\varphi_{m,j}^*)^*(-x)|^2/(2\pi)^2$ with by applying Parseval's Formula
\begin{eqnarray*}
 \sum_{j\in{\mathbb Z}}|\langle t, \varphi_{m,j}\rangle |^2\sum_{j\in{\mathbb
Z}}|(\varphi_{m,j}^*)^*(-x)|^2/(2\pi)^2
=\|t\|_2^2D_m\int
\varphi^*(u)^2du/(2\pi) = D_m\|t\|_2^2.\end{eqnarray*}

Let $b$ such that $1/2<b <a_g$.  Since $ \|g-g_m\|_{\infty} \leq
(2\pi)^{-1}\int_{|x|\geq \pi D_m} |g^*(x)|dx$
we get that \begin{eqnarray*}
\|g-g_m\|_{\infty} &\leq&
(2\pi)^{-1}((\pi D_m)^2+1)^{-(a_g-b)/2} e^{-B_g |\pi
D_m|^{r_g}} \int_{|x|\geq \pi D_m}
|g^*(x)|(x^2+1)^{(a_g-b)/2} e^{B_g |x|^{r_g}} dx\end{eqnarray*}
also bounded by
\begin{multline*}
\frac 1{2\pi}((\pi D_m)^2+1)^{-(a_g-b)/2} \exp(-B_g |\pi
D_m|^{r_g})
\kappa_{a_g}^{1/2} \sqrt{ \int_{|x|\geq \pi D_m} (x^2+1)^{-b} dx}
\\
\leq(2\pi)^{-1}((\pi D_m)^2+1)^{-(a_g-b)/2} \exp(-B_g
|\pi D_m|^{r_g}) \kappa_{a_g}^{1/2}(\pi D_m)^{1/2-b}\\
\leq  (2\pi)^{-1}\sqrt{\pi D_m} ((\pi D_m)^2+1)^{-a_g/2}
\exp(-B_g |\pi D_m|^{r_g}) \kappa_{a_g}^{1/2}.\end{multline*}
\hfill $\Box $


\end{document}